\documentclass[11pt,a4paper,intlimits]{amsart}

\usepackage{amsmath}
\usepackage{esint}
\usepackage{amsthm}
\usepackage{amssymb}
\usepackage{amsfonts}
\usepackage{amsxtra}
\usepackage{color}
\usepackage{enumitem}
\usepackage{pdfsync}

\usepackage{graphicx}
\usepackage{type1cm}
\usepackage{eso-pic}
\usepackage{color}
\makeatletter
\AddToShipoutPicture{%
            \setlength{\@tempdimb}{.05\paperwidth}%
            \setlength{\@tempdimc}{.5\paperheight}%
            \setlength{\unitlength}{1pt}%
            \put(\strip@pt\@tempdimb,\strip@pt\@tempdimc){%
        \makebox(0,0){\rotatebox{90}{\textcolor[gray]{0.75}%
        {\fontsize{1cm}{1cm}\selectfont{Draft version \quad  - \quad
 \today} }}}%
            }%
}
\numberwithin{equation}{section}

\newtheorem{theorem}{Theorem}[section]
\newtheorem{lemma}[theorem]{Lemma}
\newtheorem{proposition}[theorem]{Proposition}

\newtheorem{definition}[theorem]{Definition}

\newtheorem{remark}[theorem]{Remark}

\newcommand{\pa}{\partial}
\newcommand{\embedding}{\hookrightarrow}
\newcommand{\norm}[1]{\left\lVert #1\right\rVert}
\newcommand{\setR}{\mathbb{R}}

\DeclareMathOperator{\grad}{grad}
\DeclareMathOperator{\supp}{supp}
\DeclareMathOperator{\id}{id}
\DeclareMathOperator{\dv}{div}
\DeclareMathOperator{\loc}{loc}
\DeclareMathOperator{\sgn}{sgn}
\DeclareMathOperator{\tv}{TV}
\DeclareMathOperator{\BV}{BV}
\DeclareMathOperator{\ric}{Ric}
\DeclareMathOperator{\prob}{Prob}
\DeclareMathAlphabet{\mathpzc}{OT1}{pzc}{m}{it}

\begin{document}

\title[Scalar conservation laws on Riemannian manifolds]{Scalar conservation laws on constant and
time-dependent Riemannian manifolds}

\author{%
  {%
    Daniel~Lengeler$^*$, %
    Thomas~M\"uller$^{**}$ %
  }
   \\[1em]%
   \vspace{0.3cm}%
   \\
   {%
     \scriptsize%
     $^*$Fakult\"at f\"ur Mathematik, Universit\"at Regensburg,\\
     Universit\"atsstr. 31, D-93053 Regensburg, Germany.%
   }
   \\[1em]%
  {%
    \scriptsize%
    E\lowercase{mail}: \lowercase{\tt daniel.lengeler@mathematik.uni-regensburg.de}%
  }%
  \\
  \\
  {%
    \scriptsize%
    $^{**}$Abteilung f\"ur Angewandte Mathematik, Universit\"at Freiburg,\\
    Hermann-Herder-Str. 10, D-79104 Freiburg, Germany. %
  }\\[1em]%
  {%
    \scriptsize%
    E\lowercase{mail}: \lowercase{\tt mueller@mathematik.uni-freiburg.de}
   }\\[1em]%
}

% \author{Daniel Lengeler}
% \address{Fakult\"at f\"ur Mathematik, Universit\"at Regensburg, Universit\"atsstr. 31
% 93053 Regensburg, Germany.}
% \email{daniel.lengeler@mathematik.uni-regensburg.de}
% 
% \author{Thomas M\"uller}
% \address{Abteilung f\"ur Angewandte Mathematik, Universit\"at Freiburg,
% Hermann-Herder-Str. 10, 79104 Freiburg, Germany.}
% \email{mueller@mathematik.uni-freiburg.de}

\thanks{The work has been supported by
Deutsche Forschungsgemeinschaft via SFB/TR 71 `Geometric Partial
Differential Equations'.}

\begin{abstract}
In this paper we establish well-posedness for scalar conservation laws on closed manifolds $M$ endowed with a constant or a
time-dependent Riemannian metric for initial values in $L^\infty(M)$. In particular we show the existence and uniqueness of
entropy solutions as well as the $L^1$ contraction property and a comparison principle for these solutions.
Throughout the paper the flux function is allowed to depend on time and to have non-vanishing divergence. Furthermore, we derive estimates of the total variation of
the solution for initial values in $\BV(M)$, and we give, in the case of a time-independent metric, a simple
geometric characterisation of flux functions that give rise to total variation diminishing estimates.
\end{abstract}

\keywords{Hyperbolic; Conservation laws; Riemannian manifolds; Measure-valued solutions; Shock waves}

\subjclass{35L65; 58J45; 76N10}

\maketitle

\pagestyle{myheadings}
\thispagestyle{plain}
\markboth{Daniel Lengeler, Thomas M\"uller}{Scalar conservation laws on Riemannian manifolds}

 \section{Introduction}
% 
% \marginpar{Anwendungen, beinhaltet Hyperflaechen ,Cauchyproblem, bisherige Literatur (Abgrenzung von BenArtzi \& LeFloch), Zusammenfassung der Kapitel}
% 
Many physical systems from continuum mechanics can be modelled by nonlinear conservation laws. Typically, these partial differential equations of hyperbolic type are posed in Euclidean space. For some applications, though, the suitable domains turn out to be hypersurfaces or, more generally, Riemannian manifolds which, additionally, may change in time. Consider for example the shallow water equations on the sphere \cite {WDHJS92} as a model for the global air and water flow. Further examples include the flow of oil on a moving water surface, the transport processes on cell surfaces \cite{AS11,RS05}, surfactants on the interfacial hypersurface between two phases in multiphase flow \cite{BPS05}, the flow of a fluid in fractured porous media whose fractions are considered as a lower dimensional manifold \cite{GLSW10},
and relativistic flows \cite{Gil00,SBG01}. Scalar conservation laws have been established as a good model problem for studying the nonlinear effects in such systems.

In this contribution we consider a nonlinear scalar conservation law posed on a closed\footnote{i.e. compact, without a boundary} manifold $M$ endowed with a constant or time-dependent Riemannian metric $g$ on a given time interval $I$. The conserved quantity $u:M\times I\rightarrow\setR$ is transported according to a flux function $f$ and undergoes compression and rarefaction due to the time-dependence of $g$. Here $f=f(x,t,u)$ is a general family of vector fields on $M$ parametrized by time $t$ and by the real parameter $u$; in particular we do not assume $f$ to be divergence-free. In integral form the conservation law reads
\begin{equation}\label{eqn:conlawint}
   \frac{d}{dt} \int_K u \; dV +  \int_{\pa K} g(f(u),n_{\pa K}) \; dV_{\pa K}  = 0 \text{\quad in } I 
\end{equation}
 for some smoothly bounded $K\subset M$. Here, $dV$ and $dV_{\pa K}$ are the respective volume measures of $M$ and $\pa K$, and $n_{\pa K}$ is the outer unit normal to $\pa K$. Let the scalar field $\lambda:M\times I\rightarrow\setR$ be defined by the relation $\pa_t dV = \lambda\, dV$. Applying the Gauss theorem to \eqref{eqn:conlawint} we see that the propagation of the conserved quantity $u$ is governed by 
the Cauchy problem
\begin{equation}\label{eqn:conlaw}
 \begin{aligned}
  \pa_t u + \lambda\, u + \dv f(u)  &= 0 && \text{\quad in } M\times I,\\
u(0)&=u_0 && \text{\quad in } M.
 \end{aligned}
\end{equation}
Here, $\dv$ denotes the Riemannian divergence operator, $\lambda$ accounts for the rarefaction and compression due to the time-dependence of the metric, and $u_0:M\rightarrow \setR$ is a given initial datum. Obviously, this setting includes the case of a scalar conservation law posed on a closed, moving hypersurface in $\setR^{d+1}$. This can be seen by pulling back the flux function and the metric tensor to the initial hypersurface. We show the existence and uniqueness (in the space of measure-valued entropy solutions) of entropy solutions for initial values in $L^\infty(M)$ by the vanishing viscosity method. To this end we construct smooth solutions of parabolic regularizations of the conservation law and show uniform bounds of the $L^\infty$ norms of these solutions. Due to the time-dependence of the metric and the fact that $f$ is not divergence-free the $L^\infty$ norm might increase with time. With this uniform bound at hand, letting the regularization parameter tend to zero we obtain a measure-valued entropy solution to the conservation law. The key step in showing that the measure-valued solution is, in fact, a uniquely determined scalar-valued solution is the proof of an averaged contraction property for arbitrary measure-valued solutions, 
first introduced by DiPerna \cite{diP85}. The proof we give is based on Kruzkov's doubling of variables technique. This approach seems to be new in the context of measure-valued solutions and is based on a careful application of Lebesgue's differentiation theorem, cf. the discussion in \cite{EGH95}. From the averaged contraction property we also deduce the $L^1$ contraction property as well as a comparison principle. In the last part of our paper we derive estimates of the total variation of the entropy solution, assuming the initial datum to be in $\BV(M)$. In doing so it is crucial to work in a global, coordinate-independent way in order to get sharp estimates. On a manifold the total variation of the solution may increase with time, contrary to the situation with constant flux functions in Euclidean space. For a constant-in-time metric, though, our estimates allow us to show that a flux function gives rise to total variation diminishing estimates if and only if, essentially, it is a family of Killing 
fields. Based on the total variation estimates we establish Lipschitz continuity of the solution in time. The validity of total variation estimates is important in order to obtain a-priori error estimates for finite volume schemes, cf. \cite{LON09,GW11}. Throughout the whole paper we try to minimize the regularity assumptions on $g$ and $f$. In doing so some arguments get slightly more involved. To the best knowledge of the authors, the given proofs of the strong convergence of vanishing viscosity approximations and of the total variation estimates are the first complete ones in the context of manifolds, even in the case of a constant-in-time metric.

The first paper dealing with the existence and uniqueness for nonlinear transport equations on (time-independent) manifolds was published by Panov \cite{Pan97}. In this paper a variant of \eqref{eqn:conlaw} is considered which is independent of the geometry. Due to this independence Panov can reduce the whole problem to the Euclidean case. Recently, numerous beautiful results on the theoretical and numerical analysis of conservation on manifolds were published by Ben-Artzi, LeFloch, and their collaborators, e.g. \cite{ABL05,ALO08,BFL09,BL07,LO08.2,LON09}. In the seminal contribution \cite{BL07} a scalar conservation law on a time-independent Riemannian manifold with a divergence-free, so called geometry-compatible, flux function is considered. The authors present results on the existence and uniqueness of entropy solutions, based on the averaged contraction property, as well as total variation estimates. There are two issues in this paper which need to be clarified, though. The proof of the averaged 
contraction property 
proceeds as in \cite{diP85,Sze89} 
and implicitly assumes the existence of very special approximations of general measure-valued solutions. While the existence of such approximations is a trivial fact when dealing with constant flux functions in Euclidean space, this is not so clear in the case of manifolds. Furthermore, in the proof of the total variation estimates the authors proceed in a local, coordinate-dependent way, and a problem occurring then is the fact that the commutator of the Laplacian with a non-Killing vector field is a second order differential operator (cf. \cite{Tay96}, e.g.). 
% Amorim et al. \cite{ABL05} introduce the notion of a total variation along vector fields and give a simple condition implying that the total variation along a vector field diminishes. Here, the same problem involving the commutator occurs again. Furthermore, the convergence of a finite volume scheme is proved. 
Moreover, in \cite{Pan11} Panov proves the existence and uniqueness for the Dirichlet problem of a nonlinear transport equation on a manifold without a Riemannian metric. His approach is based on a kinetic formulation. The present paper will appear simultaneously with \cite{DKM12}. In this paper the authors show the existence of entropy solutions on moving hypersurfaces in $\setR^{d+1}$ (which is a special case of our setting) using a differential calculus which is close to numerical analysis. Furthermore, they show the uniqueness of these solutions by localizing Kruzkov's doubling of variables technique. Finally, they compute approximate solutions using a finite volume scheme.

The present paper is organized as follows. Notation and preliminary results are introduced in Section \ref{sec:not}.
In Section \ref{sec:parreg} we consider the parabolic regularization of the conservation law. Section \ref{sec:wellposedness} then deals with the existence and uniqueness of entropy solutions, while in Section \ref{sec:tvestimates} we establish the total variation estimates.

\vspace{0.5cm}
\section{Notation and preliminary results}\label{sec:not}
For a general introduction to Riemannian geometry see, e.g., \cite{doC92,Lee97}. Let $M$ be a smooth, oriented $d$-dimensional manifold. With a given Riemannian metric $g=\langle\,\cdot,\cdot\,\rangle$, assumed to be sufficiently smooth, we associate the Levi-Civita connection (or covariant
derivative)
$\nabla=\nabla^g$ on the
tangential
bundle $TM$. In fact, $g$ and $\nabla$ induce a scalar product and a connection on all tensor bundles over $M$.
For example, for covector fields $\alpha,\beta$ we have\footnote{Throughout the paper we use Einstein's
summation convention.} $\langle\alpha,\beta\rangle=g^{ij}\,\alpha_i\,\beta_j$, where $(g^{ij})$ denotes the matrix inverse of $(g_{ij})$, and $\nabla
\alpha(X,Y)=X(\alpha(Y)) - \alpha(\nabla_X Y)$ for tangential vectors $X,Y$. We denote the associated norm by
$|\alpha|=|\alpha|_g$. In coordinates we have $\alpha_{j;i}=(\nabla \alpha)_{ij}=\pa_i\alpha_j-\Gamma_{ij}^k\,\alpha_k$ with
the Christoffel symbols $\Gamma_{ij}^k$. For a scalar function $u$, i.e. a $0$-tensor field, the covariant derivative is simply
the differential. In coordinates we have $u_{;i}=(\nabla u)_i=\pa_i u$. The vector field $\grad u$ is associated with the covector
field $\nabla u$ by raising the index, in coordinates $u_{;}^{\,i}=(\grad u)^i=g^{ij}\,u_{;i}$. Furthermore, the adjoint of
$\nabla$ is the divergence operator $\dv=\dv_g$. For example, in coordinates for a vector field $X$ we have $\dv X=X^i_{\;;i}$.
By $\nabla^k$ we denote the $k$-fold application of $\nabla$ which increases the covariant order of the tensor field by $k$. $g$
and $\nabla$ induce Laplace operators
acting on arbitrary tensor fields. For example, in coordinates for a covector field $\alpha$ we have
$(\Delta_g\alpha)_k=(\Delta\alpha)_k=g^{ij}\,\alpha_{k;ij}=g^{ij}\,(\nabla^2\alpha)_{jik}$. We let $dV=dV_g$ denote the Riemannian
volume form, in positively oriented coordinates $r\in\setR^d$ we have $dV= (\det(g_{ij}))^{\frac 1 2}\, dr$. When dealing with objects depending on 
the Riemannian metric (like differential operators, the volume form, etc.) we omit the index $g$ whenever this shouldn't lead to confusion. We shall make use of Riemannian normal coordinates
centred at some point $x\in M$. These are given by concatenating the inverse of the exponential map $\exp_x=\exp_x^g$ with the isomorphism from
the tangent space $T_x M $ to 
$\setR^d$ induced by choosing a $g$-orthonormal basis of $T_x M$. In these coordinates at
$x=0$ we
have $g_{ij}=\delta_{ij}$, $\pa_k g_{ij}=0$ and $\Gamma_{ij}^k=0$. For $x,y\in M$ we denote by $d(x,y)=d_g(x,y)$ the
Riemannian distance from $x$ to $y$. We let $\gamma_{xy}$ be an arc length parametrized, minimizing geodesic
connecting $x,y$ which is uniquely determined provided that $d(x,y)$ is sufficiently small. 
Furthermore, we denote by $P_{xy}:T_x M\rightarrow T_yM$ the linear isomorphism which is given by parallel transport along $\gamma_{xy}$. 

% In local coordinates around $x\in M$, for $X\in T_xM$ we have for $g\in C^{1,1}$ and $d(x,y)$ sufficiently small
% \marginpar{So m\"usste $g\in C^{1,1}$\\ gen\"ugen, oder?}
% \begin{equation}\label{eq:ineqpartrp}
%    \sum_k|(P_{xy}X)^k - X^k + \Gamma_{ij}^k(x)\; X^i\; \dot{\gamma}^j_{xy}(0) \; d(x,y) | = \mathcal{O}(d(x,y)^2).
% \end{equation}
% %where the constant depends on the Lipschitz constants of the Christoffel symbols.
% %This estimate says that close to the origin of normal coordinates parallel transport
% %doesn't differ much from a simple translation.
% In order to prove this we let
% $\eta^k(s):=(P_{x\gamma_{xy}(s)}X)^k - X^k + \Gamma_{ij}^k(x)\; X^i\; \dot{\gamma}^j_{xy}(0) \; s$. 
% Obviously, $\eta^k(0)=0$ and in view of the differential equation defining parallel transport (see e.g. \cite{doC92}) 
% %it is easy to see that
% also $\frac{d\eta^k}{ds}(0) = 0 $. Then \eqref{eq:ineqpartrp} follows
% since $\eta$ is $C^{1,1}$.\marginpar{$\eta$ ist $C^{1,1}$, oder?}

%\begin{equation*}
% \frac{d\eta^k}{ds}(s) = -\dot{\gamma}_{xy}^i(s)\, (\eta^j(s)+X^j)\, \Gamma^k_{ij}(\gamma_{xy}(s)).
%\end{equation*}
%Then, \eqref{eq:ineqpartrp} follows from $\Gamma^k_{ij}(\gamma_{xy}(s))=\mathcal{O}(s)$, $\eta(0)=0$, and an application of the Gronwall
%inequality.

Let $U$ be a coordinate patch of $M$ with smooth coordinate vector fields $\pa_i$. Applying the Gram-Schmidt algorithm
to these fields we obtain continuous vector fields $e_i$ which form a $g$-orthonormal basis of
$T_xM$ for each $x\in U$. To be more precise, we define inductively
\[e_i:=\frac{\pa_i-\sum_{k=1}^{i-1}\langle\pa_k,e_k\rangle\,e_k}{|\pa_i-\sum_{k=1}^{i-1}\langle\pa_k,e_k\rangle\,e_k|}.\]
These fields obviously have the same regularity as $g$ and are called an orthonormal frame.

For a continuous Riemannian metric, all integrable functions $f$, and almost all $x\in M$ we have
\begin{equation*}
\lim_{\epsilon\searrow 0}\fint_{\mathcal{B}_\epsilon(x)}|f(y)-f(x)|\ dV(y)=0,
\end{equation*}
i.e. almost all $x\in M$ are Lebesgue points of $f$. Here $\mathcal{B}_\epsilon(x)$ is the geodesic ball of radius $\epsilon$
centred at $x$ and $\fint$ is the mean value integral. This fact can be easily proved in local coordinates by using the facts
that the Riemannian metric is comparable to the Euclidean metric and that in Euclidean space for a fixed integrable
function almost every point is a Lebesgue point. Now, let us state the
following simple lemma.
\begin{lemma}\label{lem:lebesgue}
Let $K$ be a smooth manifold and let $M$ be endowed with a continuous Riemannian metric. Furthermore, let $f=f(x,k): M\times
K\rightarrow
\setR$ be integrable in $x$ for all $k$ and continuous in $k$ for almost all $x$. We assume that for all $(x_0,k_0)$ there exists
a neighbourhood $U$ of $x_0$ such that the continuity in $k_0$ is uniform w.r.t.\footnote{with respect to} all $x\in U$. Then
almost all $x\in M$ are Lebesgue points for all $k\in K$. 
\end{lemma}
\proof We choose a countable, dense subset $K_0\subset K$. Obviously, the claim is true for all $k_0\in K_0$. For arbitrary
$k\in K$ we have
\begin{equation*}
 \begin{aligned}
  \fint_{\mathcal{B}_\epsilon(x)}|f(y,k)-f(x,k)|\ dV(y)\le&\fint_{\mathcal{B}_\epsilon(x)}|f(y,k_0)-f(x,k_0)|\ dV(y)\\
& +\fint_{\mathcal{B}_\epsilon(x)}|f(y,k)-f(y,k_0)|\ dV(y)\\
& + |f(x,k)-f(x,k_0)|.
 \end{aligned}
\end{equation*}
The last two terms on the right hand side get small if we choose $k_0$ close to $k$, while for fixed $k_0$ the first term tends to
zero for almost all $x$ and $\epsilon\searrow 0$.
\qed\\

We fix a function $\rho\in C^\infty_c(-1,1)$ such that $\rho\ge 0$ and $\int_{\setR^d}\rho(|x|)\ dx=1$, and for $x,y\in M$ we set
$\rho_\epsilon(x,y):=\rho_\epsilon(d(x,y)):=\frac{1}{\epsilon^d}\rho\big(\frac{d(x,y)}{\epsilon}\big)$. 
% Let us fix coordinates $(r^i)$ centred at $x$. 
% Applying the Gram-Schmidt algorithm to the standard basis of $\setR^d$ w.r.t. $g_{ij}(0)$
% and replacing the standard basis with the resulting basis we may assume without loss of generality that
% $g_{ij}(0)=\delta_{ij}$. Using these coordinates we see that for $\epsilon\searrow 0$
Let us fix coordinates $(r^i)$ centred at $x$ and apply the Gram-Schmidt algorithm to the standard basis of $\setR^d$ w.r.t. $g_{ij}(0)$. Concatenating the original chart with the isomorphism that maps the standard basis to the resulting basis we may assume without loss of generality that $g_{ij}(0)=\delta_{ij}$.\footnote{In the following we shall employ these coordinates instead of Riemannian normal coordinates in order to minimize the regularity assumptions on $g$.} In these coordinates we have
% Using these coordinates we see that for $\epsilon\searrow 0$
% \begin{equation}\label{eqn:masseeins}
%  \begin{aligned}
% \int_M\rho_\epsilon(x,y)\ dV(y)=\int_{\setR^d}
% \rho_\epsilon(d(0,r))\,(\det(g_{ij}(r)))^{\frac 1 2}\ dr\rightarrow 1
%  \end{aligned}
% \end{equation}
% since $(\det(g_{ij}(r)))^{\frac 1 2}=1+\mathcal{O}(|r|)$ and
\begin{equation}\label{eq:quasinormalcoordinates1}
 d(0,r)=|r|+\mathcal{O}(|r|^2)
\end{equation}
for $r\rightarrow 0$ and $g\in C^{0,1}$.  In order to prove \eqref{eq:quasinormalcoordinates1} let us make a general remark that is very easy to prove. If the Riemannian metrics $g_0,g_1$ are comparable with constants $c_0,c_1$, i.e. $c_0\,g_0\le g_1\le c_1\,g_0$, then the corresponding distance functions are comparable with constants $c_0^{1/2},c_1^{1/2}$. Now, in the constructed coordinates the metric $g_{ij}$ is comparable to the Euclidean metric with constants of the form $1+\mathcal{O}(|r|)$, depending on the Lipschitz constant of $g$. This proves \eqref{eq:quasinormalcoordinates1}.

Furthermore, let us show that in the constructed coordinates the minimizing geodesic connecting 
%the points
$r$ and $0$ points at the origin, up to an error of the order $\mathcal{O}(|r|)$, i.e.
\begin{equation}
 \label{eq:comparisongeodesic-r}
\dot\gamma_{r0}(0)=-\frac{r}{|r|} + \mathcal{O}(|r|),
\end{equation}
for $r\rightarrow 0$ and $g\in C^{1,1}$. To this end consider the difference $\beta$ between $\gamma_{r0}$ and the Euclidean straight line
from $r$ to $0$,
\begin{equation*}
\beta: \;s\;\mapsto\; \gamma_{r0}(s)\,-\,(d(0,r)-s)\frac{r}{d(0,r)}.
\end{equation*}
As the curve $\beta$ meets the origin at $s=0$ and $s=d(0,r)$ we know by the mean value theorem that the derivatives
of the components $\beta^i$ vanish for some ${\bar s}^i\in [0,d(0,r)]$.
Taking into account the Lipschitz continuity of $\dot\beta$ we can complete the proof of \eqref{eq:comparisongeodesic-r} by an application of \eqref{eq:quasinormalcoordinates1}.

% In order to show \eqref{eq:quasinormalcoordinates1} we consider the local representative $\gamma$ of the 
% minimizing geodesic connecting $0$ and $r$ which is supposed to be parametrized by arc length.
% Noting that $\norm{\cdot}_{g(0)}=|\cdot|$ and using the Lipschitz continuity of $g$ 
% % + O (d(0,\gamma(s)))$ 
% we see that
% \begin{equation*}
%  |r| \leq \int_0^{d(0,r)} |\dot \gamma(s) |\;ds \leq d(0,r) + \mathcal{O}(d(0,r)^2) 
% \end{equation*}
% Analogously, using the curve $s\mapsto s\frac{r}{|r|}$, we deduce $d(0,r)\leq |r| + O (|r|^2) \leq c\; |r|$
% and the combination of both inequalities implies \eqref{eq:quasinormalcoordinates1}.

Let us suppose from now on that the metric depends on time $t\in I$ for some interval $I$. If $g\in C^{0,1}$, then  $g(t_0),g(t_1)$ are comparable with constants of the form $1+\mathcal{O}(|t_0-t_1|)$, and we deduce as above that
\[d_{g(t_0)}(x,y)=d_{g(t_1)}(x,y)+\mathcal{O}(d_{g(t_1)}(x,y)\,|t_0-t_1|),\]
in particular
\begin{equation}\label{eqn:timederivdist}
 \begin{aligned}
  \pa_t d_{g(t)}(x,y)= \mathcal{O}(d_{g(t)}(x,y)).
 \end{aligned}
\end{equation}
Of course, the function $\rho_\epsilon=\rho_\epsilon(x,y,t)$ depends on time, too. Let us again fix coordinates centred at $x\in M$. Applying the Gram-Schmidt algorithm w.r.t. $g_{ij}(0,t)$ as above we obtain time-dependent coordinates such that $g_{ij}(0,t)=\delta_{ij}$. Using these coordinates we see that for $\epsilon\searrow 0$
\begin{equation}\label{eqn:timederiv}
 \begin{aligned}
  \pa_t\int_M\rho_\epsilon(x,y,t)\ dV_{g(t)}(y)=\pa_t\int_{\setR^d}
\rho_\epsilon(d(0,r;t))\,(\det(g_{ij}(r,t)))^{\frac 1 2}\ dr
\rightarrow 0,
 \end{aligned}
\end{equation}
since $\pa_t g_{ij}(0,t)=0$ and 
\begin{equation}\label{eq:quasinormalcoordinates2}
\pa_t d(0,r;t)=\mathcal{O}(|r|^2) 
\end{equation}
for $r\rightarrow 0$ and $g\in C^{0,1}$. In order to prove \eqref{eq:quasinormalcoordinates2} we just note that $g_{ij}(t_0),g_{ij}(t_1)$ are comparable with constants of the form $1+\mathcal{O}(|r|\,|t_0-t_1|)$.\footnote{It would be slightly easier to prove \eqref{eqn:timederiv} with the help of normal coordinates. But the technique employed is needed for \eqref{eqn:timederivdist} anyway and, furthermore, the use of normal coordinates would need $g$ to be more regular.}

Now, let us define function spaces on $M$ under the assumption that $M$ is closed. To this end we fix some finite atlas $(\varphi_k,U_k)$ with a subordinate partition of unity
$(\psi_k)$. We say that a measurable $(l,r)$-tensor field $T\in X(M)$ if
$T_{i_1,\ldots,i_r}^{j_1,\ldots,j_l}\,\psi_k\circ\varphi_k^{-1}\in X(\varphi_k(U))$ for
all $k,i_n,j_m$, where $X=L^p$ or $W^{k,p}$. The norm is defined by
\[\norm{T}_{X(M)}:=\sum_{k,i_n,j_m} \norm{T_{i_1,\ldots,i_r}^{j_1,\ldots,j_l}\,\psi_k\circ\varphi_k^{-1}}_{X(\varphi_k(U))}.\] 
Obviously, these spaces do not depend on the specific choice of the atlas and the partition of unity. Furthermore, the above norms
are equivalent to the canonical norms w.r.t. the metric $g=g(t)$ (uniformly w.r.t. $t$ on bounded time intervals) provided that
$g$ is sufficiently regular, e.g.
\[\norm{T}_{H^{k}(M)}^2\sim \sum_{j=0}^k\int_M |\nabla^j T|^2\ dV,\]
where $H^k:=W^{k,2}$. We similarly define the spaces $L^p(M\times I)$ and $L^p(M\times I\times\setR)$. Whenever we wish to emphasize that given function spaces are normed w.r.t. a specific metric we write $L^p(M,g)$, etc.

Finally, we provide a simple technique to produce mollifications of general tensor fields. Let $T$ be a measurable $(r,l)$-tensor field on $M$. We define the components w.r.t. the chart
$\varphi_k$ of the locally supported tensor field $T_{\epsilon,k}$ to be
\[(T_{\epsilon,k})_{i_1,\ldots,i_r}^{j_1,\ldots,j_l}:=\omega_\epsilon\ast
T_{i_1,\ldots,i_r}^{j_1,\ldots,j_l}\ \psi_k\circ\varphi_k^{-1},\]
and we set $T_\epsilon:=\sum_k T_{\epsilon,k}$. Here $\omega_\epsilon$ denotes a standard $d$-dimensional mollifier kernel.
$T_\epsilon$ is a smooth, globally defined $(r,l)$-tensor field. By Euclidean theory we have $T_\epsilon\rightarrow T$ in $L^p(M)$ if $1\le p<\infty$ and $T\in L^p(M)$, and
$T_\epsilon\rightarrow T$ uniformly if $T$ is continuous. Furthermore, for scalar fields $T\in L^\infty(M)$ we have $\norm{T_\epsilon}_{L^\infty(M)}\rightarrow\norm{T}_{L^\infty(M)}$.
This follows, like in the Euclidean case, from the weak-$^*$ lower semicontinuity of the $L^\infty$ norm and the fact that $\norm{T_\epsilon}_{L^\infty(M)}\le \norm{T}_{L^\infty(M)}$. We similarly define mollifications of tensor fields depending additionally on real parameters like time.

\vspace{0.5cm}
\section{Parabolic regularization}
\label{sec:parreg}
For the rest of the paper we fix a smooth, oriented, closed $d$-dimensional manifold M, a time interval $I=(0,T)$, a time-dependent Riemannian metric $g$ on $M$, and a vector field $f$ on $M$ depending on time and on a real parameter $u$. Throughout this section we assume $g,f$ to be smooth and that $f$ and its derivatives of any order are uniformly bounded. In particular, $f$ is globally Lipschitz continuous in $u$. Furthermore, in this section we let the function spaces be normed w.r.t. a fixed Riemannian
metric, e.g. $g(t=0)$. We intend to construct a smooth,
real-valued solution $u$ of the nonlinear parabolic equation
\begin{equation}\label{eqn:pareg}
 \begin{aligned}
  \pa_t u + \lambda\, u + \dv f(u)  &= \epsilon\,\Delta u && \text{\quad in } M\times I,\\
u(0)&=u_0 && \text{\quad in } M
 \end{aligned}
\end{equation}
for $\epsilon>0$ and $u_0\in C^\infty(M)$. Note that $\lambda$ and the differential operators involved depend on the metric and hence on time. It is not hard to check that locally $\lambda=-g^{ij}\,\pa_t g_{ij}$.

\begin{proposition}\label{prop:smoothsol}
 For every $\epsilon>0$ and every $u_0\in C^\infty(M)$ \eqref{eqn:pareg} admits a unique smooth
solution.
\end{proposition}
\proof We intend to apply the contraction mapping principle. To this end we first need to construct
appropriate solutions of the linear equation
\begin{equation*}
 \begin{aligned}
  \pa_t u + \lambda\, u + \dv f(v)  &= \epsilon\,\Delta u && \text{\quad in } M\times
I,\\
u(0)&=u_0 && \text{\quad in } M
 \end{aligned}
\end{equation*}
for given $v\in C(\bar I,L^2(M))$. To this end we set $H:=L^2(M)$, $V:=H^1(M)$, $b_v:=-\dv
f(v)\in L^2(I,V')$, and
\[a(w_0,w_1):=\int_M \epsilon\,\langle dw_0,dw_1\rangle + \lambda\,w_0\,w_1\ dV\]
for $w_0,w_1\in V$. Note that $a(\cdot,\cdot)$ depends on time since the scalar product,
 $\lambda$ and the volume form do so. Obviously, we have
\begin{equation}\label{eqn:a}
 \begin{aligned}
|a(w_0,w_1)|\le c_0\norm{w_0}_{V}\norm{w_1}_{V},\\
a(w_0,w_0)\ge c_1\norm{w_0}_V^2 - \lambda_0\norm{w_0}_H^2
 \end{aligned}
\end{equation}
for $w_0,w_1\in V$ and positive constants $c_0,c_1,\lambda_0$. By Theorem 4.1 and Remark 4.3 in
\cite[Chapter 3]{Lio72} there exists a unique $u\in L^2(I,V)$ admitting a distributional derivative $u'\in
L^2(I,V')$ such that $u(0)=u_0$\footnote{This makes sense due to the
embedding $\big\{u\in
L^2(I,V)\ |\ u'\in L^2(I,V')\big\}\embedding C(\bar I,H).$} and
\[\prec\! u'(t), w\!\succ + a(u(t),w) = \prec\! b_v(t),w\!\succ\]
for almost all $t\in I$ and all $w\in V$. Here, $\prec\!\cdot,\cdot\!\succ$ denotes
the duality pairing between $V'$ and $V$. For $v,\tilde v\in C(\bar I,H)$, the corresponding solutions $u,\tilde u$, and almost every $t\in I$ we have 
\[\prec\! u'(t)-\tilde u'(t), w\!\succ + a(u(t)-\tilde u(t),w) = \prec\! b_{v}(t)-b_{\tilde
v}(t),w\!\succ.\]
Setting $w:=u(t)-\tilde u(t)$ and applying \eqref{eqn:a}$_2$ we get
\begin{equation*}
 \begin{aligned}
\frac{1}{2}\frac{d}{dt}\norm{u(t)-\tilde u(t)}_H^2 &+ c_1\norm{u(t)-\tilde u(t)}_V^2\\ 
&\le \prec\! b_v(t)-b_{\tilde v}(t),u(t)-\tilde u(t)\!\succ + \lambda_0\norm{u(t)-\tilde
u(t)}_H^2\\
&\le c\norm{v(t)-\tilde v(t)}_{H}\norm{u(t)-\tilde u(t)}_V + \lambda_0\norm{u(t)-\tilde
u(t)}_H^2.
 \end{aligned}
\end{equation*}
The second inequality follows from the global Lipschitz continuity of $f$. Integrating from $0$ to $t\in I$ and using the fact that $u(0)=\tilde u(0)$ and Young's inequality we deduce that
\begin{equation*}
 \begin{aligned}
  \norm{u(t)-\tilde u(t)}_H^2 &\le c\int_0^t \norm{v(s)-\tilde v(s)}_H^2\ ds + 2\,\lambda_0\int_0^t
\norm{u(s)-\tilde u(s)}_H^2\ ds\\
&\le c\,T\norm{v-\tilde v}^2_{C(\bar I,H)} + 2\,\lambda_0\int_0^t
\norm{u(s)-\tilde u(s)}_H^2\ ds.
 \end{aligned}
\end{equation*}
By Gronwall's lemma we have
\[\norm{u-\tilde u}_{C(\bar I,H)}\le ( c\,T (1+2\lambda_0 T \,e^{2\lambda_0 T}) )^{\frac{1}{2}}  \norm{v-\tilde v}_{C(\bar
I,H)}.\]
Replacing $I$ by a smaller interval if necessary we infer that the mapping $v\mapsto u$ from $C(\bar I,H)$ to itself is Lipschitz continuous with a Lipschitz constant smaller than $1$, hence admitting a unique fixed point. Since the size of the time interval is independent of $\norm{u_0}_H$ we can repeat this argument a finite number of times to obtain a unique solution $u$ on the whole time interval $I$. 

It remains to show that $u$ is smooth. By a
localisation argument this follows from Euclidean linear parabolic theory similarly to the proof of
Lemma 2.16 in \cite[Chapter 2]{Mal96}. Let us explain the localisation argument by carrying out the first step.
Let $\varphi\in C^\infty_c(U)$ for some coordinate patch $U\subset M$. Then the localisation
$\varphi\,u$ solves
\begin{equation}\label{eqn:local}
 \begin{aligned}
\prec\! (\varphi\,u)', w\!\succ &+ a(\varphi\,u,w)\\
& = \prec\!
b_u,\varphi\,w\!\succ-\epsilon\int_U  
u\,w\,\Delta\varphi +2w\,\langle du,d\varphi\rangle\ dV
 \end{aligned}
\end{equation}
almost everywhere in $I$ and for all $w\in V$. Since $u\in L^2(I,V)$ the functional defined by the right
hand side is in fact induced by some element of $L^2(I\times U)$. Thus, writing \eqref{eqn:local} in
local coordinates we see that $\varphi\, u$ solves a uniformly parabolic equation in some smooth
domain of $\setR^d$ with vanishing boundary values, smooth initial datum and right hand side in
$L^2$. By Theorem 6.2 in \cite[Chapter 4]{Lio72II} $\varphi\,u$ lies in $L^2(I,H^2(U))$ admitting a time derivative
in $L^2(U\times I)$. As $U$ and $\varphi$ are arbitrary we see that $u$ lies in $L^2(I,H^2(M))$
admitting a time derivative in $L^2(M\times I)$. Using this localisation argument and
differentiating the localized equation repeatedly w.r.t. space and time we can adapt the proof of
Lemma 2.16 in \cite[Chapter 2]{Mal96} to show that $u$ is
smooth. Note that we need to invoke parabolic $L^p$ theory for time-dependent operators, see, e.g., \cite{LSU68}.
\qed\\

Adapting a standard argument it is not hard to show that for smooth solutions $u,\tilde
u$ of \eqref{eqn:pareg} and all $0\le s\le t\le T$ we have
\begin{equation}\label{eqn:fundest}
 \begin{aligned}
  \int_M [u(t)-\tilde u(t)]^+\ dV_{g(t)}\le \int_M [u(s)-\tilde u(s)]^+\ dV_{g(s)}
 \end{aligned}
\end{equation}
where $[v]^+:=\max(v,0)$ for $v\in\setR$. From this inequality we can deduce an $L^1$ contraction
property as well as
\[u\le \tilde u \text{ in } \bar I\times M\text{ provided that }u(0)\le \tilde u(0) \text{ in }M.\]
Note, however, that from this comparison principle we cannot infer
uniform bounds for the solutions of \eqref{eqn:pareg} because constants don't solve this equation.
This is due to the presence of the factor $\lambda$ and the nonvanishing divergence of $f$. In order
to derive uniform bounds we need to adjust the proof of \eqref{eqn:fundest} appropriately. 

\begin{proposition}\label{prop:uniformbound}
For every smooth solution $u$ of \eqref{eqn:pareg} and all $t\in\bar I$ we have
% \footnote{For $c_2=0$
% the right hand side of \eqref{eqn:uniformbound} equals $\norm{u_0}_{L^\infty(M)}+c_3\,t$.}
\begin{equation}\label{eqn:uniformbound}
 \begin{aligned}
% \norm{u(t)}_{L^\infty(M)}\le\norm{u_0}_{L^\infty(M)} \,e^{c_2
% t}+c_3/c_2\,(e^{c_2 t}-1).\\
\norm{u(t)}_{L^\infty(M)}\le\norm{u_0}_{L^\infty(M)} \,e^{\int_0^t c_2(s)\, ds
}+\int_0^t e^{\int_s^t c_2(\tau)\, d\tau}\,c_3(s)\ ds,
 \end{aligned}
\end{equation}
where
\begin{equation*}
 \begin{aligned}
c_2(t)&:=-\inf_{x\in M}\lambda(x,t) -\inf_{(x,\bar u)\in M\times\setR} \pa_u\dv f(x,t,\bar u),\\
 c_3(t)&:=\norm{\dv f(\,\cdot\,,t,0)}_{L^\infty(M)}.
\end{aligned}
\end{equation*}
\end{proposition}
\proof We set 
\[\eta_\delta(v):=\left\{\begin{array}{cl} 0 &\text{ for
}v\le 0\\
v^2/4\delta & \text{ for } 0<v\le2\delta\\
v-\delta & \text{ for } 2\delta<v\end{array}\right..\]
Note that $\eta_\delta\in W^{2,\infty}_{\loc}(\setR)$, $\eta_\delta'\ge 0$, $\eta_\delta''\ge 0$,
$\delta\,\eta_\delta''\le 1/2$, and for $\delta\searrow 0$ and all $v\in\setR$ 
\[\eta_\delta(v)\rightarrow [v]^+,\ \eta_\delta'(v)\,v-\eta_\delta(v)\rightarrow 0,\
\eta_\delta''(v)\rightarrow 0.\]
We denote the right hand side of \eqref{eqn:uniformbound} by $\tilde u$. This function solves
\eqref{eqn:pareg} with $\tilde u(0)=\norm{u_0}_{L^\infty(M)}$, $\lambda$ replaced by $-c_2$ and
$\dv f(\tilde u)$ replaced by $-c_3$. Hence we have
\begin{equation*}
 \begin{aligned}
\pa_t\eta_\delta(u&-\tilde u) + \eta_\delta'(u-\tilde u)\, (u-\tilde u)\,\lambda +
\dv\big(\eta_\delta'(u-\tilde u)(f(u)-f(\tilde u))\big)\\
&\hspace{3.82cm}- \eta_\delta''(u-\tilde u)\, d(u-\tilde u)\big(f(u)-f(\tilde u)\big)&&\\
&\hspace{0.04cm}=\pa_t\eta_\delta(u-\tilde u) + \eta_\delta'(u-\tilde u)\,\big( (u-\tilde
u)\,\lambda +
\dv(f(u)-f(\tilde u))\big)\\ 
&\hspace{0.04cm}\le\pa_t\eta_\delta(u-\tilde u) + \eta_\delta'(u-\tilde u)\,\big(u\,\lambda+ \dv
f(u) + \tilde u\,c_2 +c_3\big)\\
&\hspace{0.04cm}=\epsilon\,\Delta\eta_\delta(u-\tilde u) - \eta_\delta''(u-\tilde u)\,|d(u-\tilde
u)|^2_g\le\epsilon\,\Delta\eta_\delta(u-\tilde u).
 \end{aligned}
\end{equation*}
Integrating this inequality over $(0,t)\times M$ and integrating by parts gives
\begin{equation*}
 \begin{aligned}
\int_M\eta_\delta(u(t)-\tilde u(t))&\ dV_{g(t)}-\int_M\eta_\delta(u(0)-\tilde u(0))\ dV_{g(0)}\\
&\le\int_I\int_M \Big(\eta_\delta''(u-\tilde u)\, d(u-\tilde u)\big(f(u)-f(\tilde u)\big)\\
&\hspace{1.5cm}+\lambda\,\big(\eta_\delta'(u-\tilde u)\, (u-\tilde u) -
\eta_\delta(u-\tilde u)\big)\Big)\, dV_{g(s)}\,ds.
 \end{aligned}
\end{equation*}
Since all integrands are uniformly bounded and converge pointwise by the dominated convergence
theorem we deduce 
\begin{equation*}
 \begin{aligned}
  \int_M[u(t)-\tilde u(t)]^+\ dV_{g(t)}\le 0.
 \end{aligned}
\end{equation*}
Hence $u\le\tilde u$. Replacing $\norm{u_0}_{L^\infty(M)}$ by $-\norm{u_0}_{L^\infty(M)}$, $c_3$ by $-c_3$ and $u-\tilde u$ by
$\tilde u- u$ we similarly show that $-\tilde u\le u$.
\qed
%On the other hand, this expression dominates the left hand side of \eqref{eqn:keyid}. Now, we can
%repeat the argument in the proof of Proposition \ref{prop:fundest} to show $u\le\tilde u$.
%Replacing $\tilde u_0$ by $-\tilde u_0$, $\nu$ by $-\nu$ and $u-\tilde u$ by $\tilde u- u$ we can
%similarly show $-\tilde u\le u$.
%\qed\

\vspace{0.5cm}
\section{Well-posedness}
\label{sec:wellposedness}
Throughout this section we assume that $g\in C^{1,1}$ and $f,\pa_u f\in C^{1}$.\footnote{It might be possible to lower the regularity assumptions. Note, however, that it is only for $g\in C^{1,1}$ that the Christoffel symbols are Lipschitz continuous and geodesics are uniquely defined. Furthermore, ODE theory tells us that for $g\in C^{1,1}$ and $x\in M$ the exponential map $\exp_x$ is a local $C^1$-diffeomorphism.} Furthermore, we assume $\pa_u \dv  f$ to be uniformly bounded from below in $M\times I\times\setR$. We will establish well-posedness of the scalar conservation law \eqref{eqn:conlaw} for initial values in $L^\infty(M)$. For arbitrary convex functions $\eta:\setR\rightarrow\setR$ we set
\begin{equation*}
q(x,t,u):=\int_k^u\eta'(v)\,\pa_v f(x,t,v)\ dv
\end{equation*}
for some $k\in\setR$. For $M(\setR):=(C_0(\setR))'$ we let $\prob(\setR)\subset M(\setR)$ denote
the subset of probability measures. For a compactly supported $\nu\in \prob(\setR)$ and a continuous function $\eta$ we set
$\langle\nu,\eta\rangle:=\int_{\setR}\eta\ d\nu.$
% $\eta(u;k):=|u-k|$, $q(u;k):=\sgn(u-k)(f(u)-f(k))$.

\begin{definition}\label{def:mvsol}
For $u_0\in L^\infty(M)$ we say that a weak-$^*$ measurable\,\footnote{More precisely, for each $h\in L^1(M\times I,C_0(\setR))$ the mapping
$(x,t)\mapsto\langle\nu_{x,t},h(x,t)\rangle$ is measurable.} mapping $\nu:M\times I\rightarrow
\prob(\setR)\subset M(\setR)$, with supports $\supp\nu_{x,t}$ contained in a fixed, bounded interval, is a measure-valued
entropy solution of \eqref{eqn:conlaw} if \footnote{$\dv^x$ denotes the divergence w.r.t. to the explicit dependence on the 
variable $x$, i.e. fixing the value of $u$.}
\begin{equation}\label{eqn:mvsol}
 \begin{aligned}
\int_I\int_M \pa_t\varphi\,\langle\nu,\eta\rangle+ \varphi\,
\langle\nu,\lambda\,(\eta-\eta'(\,\cdot\,)\,\cdot\,) + \dv^x q -\eta'\dv^x f\rangle\\
 + \langle\nu,d\varphi(q)\rangle\ dV\,dt + \int_M\varphi(0)\,\eta(u_0)\ dV_{g(0)}\ge 0
\end{aligned}
\end{equation}
holds for every convex function $\eta$ and every nonnegative test function $\varphi\in
C^{0,1}_c(M\times [0,T))$. We say that $u\in L^\infty(M\times I)$ is an entropy solution if the
family of Dirac measures $\delta_u$ is a measure-valued entropy solution.
\end{definition}
\vspace{1mm}

The key step in the proof of well-posedness is to show the following averaged contraction property for measure-valued
entropy solutions. In order to get clean calculations and a clean result it is crucial to work in a global,
coordinate-independent way.
\begin{proposition}\label{prop:acp}
For measure-valued entropy solutions $\nu,\sigma$ and almost all $0<s\le t<T$ we have
\begin{equation*}
 \begin{aligned}
\int_M \langle\nu_t\otimes\sigma_t,\eta\rangle\,dV_{g(t)} \le \int_M
\langle\nu_s\otimes\sigma_s,\eta\rangle\,dV_{g(s)},
 \end{aligned}
\end{equation*}
where $\eta(u,v):=[u-v]^+$. By interchanging the roles of $\nu$ and $\sigma$ we see that the same
inequality holds for $\eta(u,v):=|u-v|$.
\end{proposition}
\proof In \eqref{eqn:mvsol} we choose $\eta=\eta(u,k):=[u-k]^+$,
\[q=q(x,t,u,k):=\sgn[u-k]^+\,(f(x,t,u)-f(x,t,k)),\]
and
\[\varphi=\varphi(x,t,y,s):=\psi(t)\,\omega_{\tilde\epsilon}(t-s)\,\rho_\epsilon(x,y,s)\]
for fixed $k\in\setR$, $y\in M$, $s\in I$. Here $\omega_{\tilde\epsilon}$ denotes a one-dimensional mollifier
kernel, $\rho_\epsilon$ is the kernel constructed in Section \ref{sec:not}, and $\psi\in C^{1}_c(0,T)$. 
% $\omega_{\tilde\epsilon}(t)=\tilde\epsilon^{-1}\omega(t\,\tilde\epsilon^{-1})$ for $\omega\in C^\infty_c(-1,1)$
% with $\int_{\setR}\omega\,dt=1$ and  
We integrate the resulting inequality over $\setR\times M\times I$ w.r.t.
$d\sigma_{y,s}(k)\,dV_{g(s)}(y)\,ds$. Interchanging the arguments of $\eta,q$ and the roles of $x,t,\nu$ and $y,s,\sigma$ we
obtain an analogous inequality. Adding these two inequalities results in
\begin{equation*}
\begin{aligned}
\int_I\int_M \psi'(t)\, I_1^{\epsilon,\tilde\epsilon}(x,t) + \psi(t)\, I_2^{\epsilon,\tilde\epsilon}(x,t) +\psi(t)\, I_3^{\epsilon,\tilde\epsilon}(x,t)\
dV_{g(t)}(x)\,dt\ge  0,
\end{aligned}
\end{equation*}
where
\begin{equation*}
\begin{aligned}
&I_1^{\epsilon,\tilde\epsilon}(x,t):=\int_I\int_M\,\omega_{\tilde\epsilon}(t-s)\,\rho_\epsilon(x,y,s)\,\langle\nu_{x
,t} \otimes\sigma_{y,s},\eta\rangle\ dV_{g(s)}(y)\,ds\\
&I_2^{\epsilon,\tilde\epsilon}(x,t):=\int_I\int_M\,\omega_{\tilde\epsilon}(t-s)\,\Big(\pa_s\rho_\epsilon(x,y,s)\,
\langle\nu_{x,t} \otimes\sigma_{y,s},\eta\rangle\\
&\hspace{2cm} + \rho_\epsilon(x,y,s)\,\big\langle\nu_{x,t}\otimes\sigma_{y,s},\lambda(x,t)\,\big(\eta(k_1,k_2) -
\pa_1\eta(k_1,k_2)\,k_1\big)\\
&\hspace{3.8cm}+\lambda(y,s)\,\big(\eta(k_1,k_2) - \pa_2\eta(k_1,k_2)\,k_2\big)\big\rangle\Big)\
dV_{g(s)}(y)\,ds\\
&I_3^{\epsilon,\tilde\epsilon}(x,t):=\int_I\int_M\omega_{\tilde\epsilon}(t-s)\,\big\langle\nu_{x,t}\otimes\sigma_{y
,t},d^x\rho_\epsilon(x,y,s)(q(x,t))\\
&\hspace{7.2cm}+d^y \rho_\epsilon(x,y,s)(q(y,s))\\
&\hspace{0.5cm}+ \rho_\epsilon(x,y,s)\,\sgn[k_1-k_2]^+\,(\dv f(y,s,k_1)-\dv f(x,t,k_2))\big\rangle\ dV_{g(s)}(y)\,ds.
\end{aligned}
\end{equation*}
Here $k_1,k_2$ denote the integration variables w.r.t. $\nu,\sigma$. Let us show that the function
$I_1^{\epsilon,\tilde\epsilon}$ tends to $I_1^\epsilon$ almost everywhere for $\tilde\epsilon\searrow 0$, where
\[I_1^\epsilon(x,t):=\int_M\rho_\epsilon(x,y,t)\,\langle\nu_{x
,t}\otimes\sigma_{y,t},\eta\rangle\ dV_{g(t)}(y).\]
To this end we note that
\begin{equation*}
 \begin{aligned}
|I_1^{\epsilon,\tilde\epsilon}(x,t)-I_1^{\epsilon}(x,t)|&=\big|\int_{\setR}\int_I\omega_{\tilde\epsilon}(t-s)\,(h(x,s,k_1)-h(x
,t,k_1))\ ds\,d\nu_{x,t}\big|\\
&\le c\,\int_{\setR}\fint_{t-\tilde\epsilon}^{t+\tilde\epsilon}|h(x,s,k_1)-h(x,t,k_1)|\ ds\,d\nu_{x,t},
 \end{aligned}
\end{equation*}
where 
\[h(x,s,k_1):=\int_M\rho_\epsilon(x,y,s)\,\langle\sigma_{y,s},\eta(k_1,k_2)\rangle\ dV_{g(s)}(y).\]
Since $h$ is continuous in $k_1$ (and $x$), uniformly w.r.t. $s$, by Lemma \ref{lem:lebesgue} the mean value integral tends to
zero for all $k_1$ and almost all $x,t$. An application of the dominated convergence theorem then proves the claim. We shall use
this \emph{convergence trick} several times throughout the proof. Using this trick
we can similarly show that $I_1^\epsilon$ tends to
\[%\psi'(t)\,
\langle\nu_{x,t}\otimes\sigma_{x,t},\eta\rangle\]
almost everywhere for $\epsilon\searrow 0$. Now, it suffices to prove that $I_2^{\epsilon,\tilde\epsilon}$ and
$I_3^{\epsilon,\tilde\epsilon}$ converge to zero almost everywhere if we let
$\tilde\epsilon\searrow 0$ and then $\epsilon\searrow 0$. 

Let us proceed with $I_2^{\epsilon,\tilde\epsilon}(t,x)$. Noting that $\eta(k_1,k_2)=\pa_1\eta(k_1,k_2)\,k_1 +
\pa_2\eta(k_1,k_2)\,k_2$ and using the Lipschitz continuity of $\lambda$ in space and time we can replace the term
$\lambda(x,t)$ by $\lambda(y,s)$ to obtain
\begin{equation*}
 \begin{aligned}
I_2^{\epsilon,\tilde\epsilon}(x,t)&=\mathcal{O}(\epsilon+\tilde\epsilon)+\int_I\int_M\omega_{\tilde\epsilon}(t-s)\,
\langle\nu_{x,t}\otimes\sigma_{y,s},\eta\rangle\\
&\hspace{3cm}\big(\pa_s\rho_\epsilon(x,y,s) + \rho_\epsilon(x,y,s)\,\lambda(y,s)\big)\
dV_{g(s)}(y)\,ds.
 \end{aligned}
\end{equation*}
By the convergence trick $I_2^{\epsilon,\tilde\epsilon}$ tends to $I_2^{\epsilon}$ almost everywhere for $\tilde\epsilon\searrow
0$, where
\begin{equation*}
 \begin{aligned}
I_2^{\epsilon}(x,t)=\int_M\langle\nu_{x,t}\otimes\sigma_{y,t},\eta\rangle\,\big(\pa_t\rho_\epsilon(x,y,t)
+ \rho_\epsilon(x,y,t)\,\lambda(y,t)\big)\ dV_{g(t)}(y)\\
+\mathcal{O}(\epsilon).
 \end{aligned}
\end{equation*}
Since
\[\big(\pa_t\rho_\epsilon(x,y,t) + \rho_\epsilon(x,y,t)\,\lambda(y,t)\big)\, dV_{g(t)}(y)=\pa_t\big(\rho_\epsilon(x,y,t)\,dV_{g(t)}(y)\big)\]
and due to \eqref{eqn:timederivdist} and \eqref{eqn:timederiv} we have
\begin{equation*}
 \begin{aligned}
|I_2^{\epsilon}(x,t)|&=\int_M\big(\langle\nu_{x,t}\otimes\sigma_{y,t},\eta\rangle-\langle\nu_{x,t}\otimes\sigma_{x,t},
\eta\rangle\big)\\
&\hspace{2cm}\big(\pa_t\rho_\epsilon(x,y,t)+ \rho_\epsilon(x,y,t)\,\lambda(y,t)\big)\ dV_{g(t)}(y) + \mathcal{O}(\epsilon)\\
&\le c\,\fint_{\mathcal{B}_\epsilon(x)}
\big|\langle\nu_{x,t}\otimes\sigma_{y,t},\eta\rangle-\langle\nu_{x,t}\otimes\sigma_{x,t},
\eta\rangle\big|\ dV_{g(t)}(y) + \mathcal{O}(\epsilon).
 \end{aligned}
\end{equation*}
Another application of the convergence trick shows that the right hand side vanishes in the limit $\epsilon\searrow 0$.

It remains to show the convergence of $I_3^{\epsilon,\tilde\epsilon}$. Using the continuity of $\dv f$ in space and time we can replace $\dv f(x,t,k_2)$ by $\dv f(y,s,k_2)$, producing a term of the order
$o(1)$ for $\epsilon+\tilde\epsilon\searrow 0$. Furthermore, note that the function
\begin{equation*}
 \begin{aligned}
&(x,s,k_1)\mapsto\int_M\rho_\epsilon(x,y,s)\,\big\langle\sigma_{y,s},\\
&\hspace{3cm}\sgn[k_1-k_2]^+\,(\dv f(y,s,k_1)-\dv f(y,s,k_2))\big\rangle\ dV_{g(s)}(y)  
 \end{aligned}
\end{equation*}
is continuous in $k_1$ (and $x$), uniformly w.r.t. $s$. Thus we can apply the convergence trick to show that
$I_3^{\epsilon,\tilde\epsilon}$ tends to $I_3^{\epsilon}$ almost
everywhere for $\tilde\epsilon\searrow 0$, where
\begin{equation*}
 \begin{aligned}
&I_3^{\epsilon}(x,t):=\int_M\big\langle\nu_{x,t}\otimes\sigma_{y,t},d^x\rho_\epsilon(x,y,t)(q(x,t))+d^y
\rho_\epsilon(x,y,t)(q(y,t))\\
&+ \rho_\epsilon(x,y,t)\,\sgn[k_1-k_2]^+\,(\dv f(y,t,k_1)-\dv f(y,t,k_2))\big\rangle\ dV_{g(s)}(y)+o(1)\\
&=:I_{3,1}^{\epsilon}(x,t)+I_{3,2}^{\epsilon}(x,t)+o(1).
 \end{aligned}
\end{equation*}
Here $I_{3,1}^{\epsilon}$ and $I_{3,2}^{\epsilon}$ comprise the terms involving $d\rho_\epsilon$ respectively $\rho_\epsilon$.
Again by the convergence trick we see that $I_{3,2}^{\epsilon}$ tends to
$I_{3,2}$ almost everywhere for $\epsilon\searrow 0$, where
\begin{equation}\label{eqn:i32}
I_{3,2}(x,t)=\big\langle\nu_{x,t}\otimes\sigma_{x,t},\sgn[k_1-k_2]^+\,(\dv f(x,t,k_1)-\dv
f(x,t,k_2))\big\rangle.
\end{equation}
For the rest of the proof we shall omit the variable $t$. Since $\grad^x d(x,y)=-\dot\gamma_{xy}(0)$ and $\grad^y d(x,y)=\dot\gamma_{xy}(d(x,y))$, provided that $d(x,y)$ is small enough, by the invariance of the Riemannian metric under parallel transport we have
\begin{equation*}
 d^x \rho_\epsilon(x,y)(q(x)) = - d^y \rho_\epsilon(x,y)(P_{xy}q(x)).
\end{equation*}
Furthermore, it's a well known fact from Riemannian geometry that
\[f(y,k_l)-P_{xy}f(x,k_l)=-\nabla_{\dot\gamma_{yx}(0)}f(y,k_l)\,d(x,y)+o(d(x,y)).\]
Thus
\begin{equation*}
\begin{aligned}
&I_{3,1}^\epsilon(x)=-\int_M\rho_\epsilon'(d(x,y))\,d(x,y)\,\big\langle\nu_{x}\otimes\sigma_{y},\sgn[k_1-k_2]^+\\
&\hspace{3cm}\langle\dot\gamma_{yx}(0),\nabla_{\dot\gamma_{yx}(0)}(f(y,k_1)-f(y,k_2))\rangle\big\rangle\ dV_{g}(y) + o(1).
\end{aligned}
\end{equation*}
Replacing $\sigma_y$ with $\sigma_x$ produces a term of the order $o(1)$. Indeed, the error term may be estimated by
\begin{equation*}
\begin{aligned}
 &c\,\fint_{\mathcal B_\epsilon(x)}|\langle\nu_x\otimes\sigma_y,\tilde h(y,k_1,k_2)\rangle - \langle\nu_x\otimes\sigma_x,\tilde h(y,k_1,k_2)\rangle | \ dV_g(y)\\
&\le c\,\fint_{\mathcal B_\epsilon(x)}|\langle\nu_x\otimes\sigma_y,\tilde h(x,k_1,k_2)\rangle - \langle\nu_x\otimes\sigma_x,\tilde h(x,k_1,k_2)\rangle | \ dV_g(y) +o(1),
\end{aligned}
\end{equation*}
where
\begin{equation*}
 %\tilde h(y,z,k_1):= \langle \sigma_z, \sgn[k_1-k_2]^+|\nabla(f(y,k_1)-f(y,k_2))|\rangle
\tilde h(y,k_1,k_2):= \sgn[k_1-k_2]^+\,\nabla(f(y,k_1)-f(y,k_2)).
\end{equation*}
An application of the convergence trick then proves the claim.
% {\bf
% Replacing $\sigma_y$ with $\sigma_x$ we produce a term that can be estimated by
% \begin{equation*}
%  R_\epsilon(x):=C \ \int_\setR \fint_{\mathcal B_\epsilon(x)}|\tilde h (y,x,k_1) - \tilde h(y,y,k_1) | \ dV_g(y) \ d\nu_x(k_1),
% \end{equation*}
% for some constant $C$ that does not depend on $\epsilon,x$ and where
% \begin{equation*}
%  %\tilde h(y,z,k_1):= \langle \sigma_z, \sgn[k_1-k_2]^+|\nabla(f(y,k_1)-f(y,k_2))|\rangle
% \tilde h(y,z,k_1):= \int_\setR \sgn[k_1-k_2]^+(\nabla(f(y,k_1)-f(y,k_2))) \ d\sigma_z(k_2).
% \end{equation*}
% Here, we used the fact that $\langle X,\nabla_X f \rangle \leq |\nabla f|$ for a vector $X$ with $|X|=1$ and a vector field $f$.
% Since $\tilde h$ is continuous in $y$ and $k_1$, uniformly w.r.t. $z$, an application of the convergence trick implies
% that $R_\epsilon(x)$ tends to zero for $\epsilon \searrow 0$ for almost all $x$.}
%use Lemma \ref{lem:lebesgue} in order to show that the mean value integral tends to zero.
Hence, in view of \eqref{eqn:i32} it suffices to show that
\[\int_M\rho_\epsilon'(d(x,y))\,d(x,y)\,\langle\dot\gamma_{yx}(0),\nabla_{\dot\gamma_{yx}(0)}f(y,k_l)\rangle\ dV_{g}(y)\ \rightarrow \ \dv f(x,k_l)\]
for $\epsilon\searrow 0$. Writing the integral in coordinates centred at $x$ with $g_{ij}(0)=\delta_{ij}$ and
taking into account \eqref{eq:quasinormalcoordinates1}, \eqref{eq:comparisongeodesic-r}, and the continuity of $\nabla f$
%, and the identity $\dot\gamma_{yx}(0)=-\frac{r}{|r|} + \mathcal{O}(|r|)$
results in
\begin{equation*}
\begin{aligned}
\sum_{i,j}f^i_{\ ;j}(0,k_l)\int_{\setR^d} \rho_\epsilon'(|r|)\,|r|\,\frac{r^ir^j}{|r|^2}\ dr + o(1).
\end{aligned}
\end{equation*}
Using polar coordinates it is easy to check that the integral equals $\delta^{ij}$. This completes the proof.\footnote{It would be slightly easier to prove this last step using normal coordinates. But this would need $g$ to be more regular.}
\qed\\

Note that the treatment of $I_{3,1}^\epsilon$ is the only step in the proof where (so far) $g$ and $f$ need to be more regular than merely Lipschitz continuous.

\begin{lemma}\label{lem:timereg}
For a measure-valued entropy solution $\nu$ we have\footnote{Here and in the following proof one has
to take $t$ (resp. $\epsilon$) from the complement of some set of vanishing measure.}
\begin{equation*}
 \begin{aligned}
 \lim_{t\searrow 0}\int_M \langle\nu_{t},|\,\cdot\,-u_0|^2\rangle\ dV_{g(t)}=0.
\end{aligned}
\end{equation*}
\end{lemma}
\proof We set
\[\varphi_{\epsilon,\tilde\epsilon}(t):=\left\{\begin{array}{cl} 1 &\text{ for
}0\le t\le \epsilon\\
1-(t-\epsilon)/\tilde\epsilon & \text{ for } \epsilon< t\le \epsilon+\tilde\epsilon\\
0 & \text{ else }\end{array}\right..\]
Choosing $\eta=\id$ and $\varphi=\varphi_{\epsilon,\tilde\epsilon}(t)\,\psi(x)$ with an arbitrary
$\psi\in C^{0,1}(M)$ in \eqref{eqn:mvsol} and letting $\tilde\epsilon\searrow 0$ we obtain for
almost every $\epsilon\in I$
\begin{equation*}
 \begin{aligned}
-\int_M \langle\nu_\epsilon,\id\rangle\,\psi\ dV_{g(\epsilon)} + \int_M u_0\,\psi\ dV_{g(0)}=
\mathcal{O}(\epsilon).
\end{aligned}
\end{equation*}
By density the left hand side vanishes in the limit $\epsilon\searrow 0$ for every $\psi\in L^1(M)$.
For $\eta(u)=u^2$ we similarly obtain
\begin{equation*}
 \begin{aligned}
-\int_M \langle\nu_\epsilon,\eta\rangle\ dV_{g(\epsilon)} + \int_M\eta(u_0)\ dV_{g(0)}\ge
\mathcal{O}(\epsilon).
\end{aligned}
\end{equation*}
These two facts imply that
\begin{equation*}
 \begin{aligned}
\limsup_{\epsilon\searrow 0}\int_M
\langle\nu_\epsilon,\underbrace{\eta-\eta(u_0)-\eta'(u_0)(\,\cdot\,-u_0)}_{=|\,\cdot\,-u_0|^2}
\rangle\ dV_{g(\epsilon)}\le 0.
%\le \limsup_{\epsilon\searrow 0} \int_M -\eta'(u_0)\langle\nu_\epsilon,\,\cdot\,-u_0\rangle\       
% dV_{g(\epsilon)}=0.
\end{aligned}
\end{equation*}
\qed\\

% \begin{remark}\label{rem:normal}
%  Let $u\in L^\infty(I\times M)$ be an entropy solution. Then the argument leading to
% \eqref{eqn:normal} shows that we may modify the function $t\mapsto u(t)$ on a set of vanishing
% measure such that for all $t\in \bar I$, all $\psi\in L^1(M)$ and $\epsilon\searrow 0$
% \[\frac{1}{\epsilon}\int_t^{t+\epsilon}\int_M u\,\psi\ dV_g\, ds \rightarrow \int_M u(t)\,\psi\
% dV_{g(t)}.\]
% From now on we assume entropy solutions to be normalized in this sense.
% \end{remark}
% \vspace{1mm}

\begin{theorem}
 For every $u_0\in L^\infty(M)$ there exists a unique entropy solution $u\in L^\infty(M\times I)$ of \eqref{eqn:conlaw}. Furthermore, $u(t)$ admits the bound \eqref{eqn:uniformbound} for almost every $t\in I$.
% \begin{equation}\label{eqn:uniformbound2}
%  \begin{aligned}
% \norm{u}_{L^\infty(I\times M)}\le\norm{u_0}_{L^\infty(M)} \,e^{c_2
% T}+c_3/c_2\,(e^{c_2 T}-1)
%  \end{aligned}
% \end{equation}
% with the constants $c_2,c_3$ defined in the last Section in terms of $g,f$.
\end{theorem}
\proof Using the procedure introduced in Section \ref{sec:not} we define mollifications $u_0^\epsilon,g^\epsilon,f^\epsilon$ of $u_0,g,f$. Furthermore, we set $f_\tau^\epsilon(x,t,u):=f^\epsilon(x,t,\tau(u))$, where the smooth function $\tau:\setR\rightarrow\setR$ with $|\tau'|\le 1$ is the identity for $|u|\le N$ and constant for large values of $u$. Here we choose $N$ to be, say, twice the maximum of the right hand side of $\eqref{eqn:uniformbound}$. From Propositions \ref{prop:smoothsol} and \ref{prop:uniformbound} we know that there exist smooth solutions $u^\epsilon$ of \eqref{eqn:pareg} admitting the bound \eqref{eqn:uniformbound} with $u_0,g,f$ replaced by $u_0^\epsilon,g^\epsilon,f^\epsilon_\tau$. Since this gives a uniform bound there exists a subsequence, again denoted by
$u^\epsilon$, converging weakly-$^*$ in $L^\infty(M\times I)$ and an associated family $\nu$ of Young measures.\footnote{For the construction of the
family of Young measures we may proceed as in the Euclidean case, see, e.g., \cite{Tar83}.} We let
$\eta:\setR\rightarrow\setR$
be a smooth, convex function. Note that for the associated entropy flux $q^\epsilon_\tau$ we have 
\begin{equation*}
\dv_{g^\epsilon} q^\epsilon_\tau(u^\epsilon)=\int_k^{u^\epsilon}\eta'(s)\,\dv_{g^\epsilon}^x\pa_u f^\epsilon_\tau(s)\
ds + \eta'(u^\epsilon)\, du^\epsilon(\pa_u f^\epsilon_\tau(u^\epsilon)).
\end{equation*}
Multiplying the identity 
\[\pa_t\eta(u^\epsilon)+\eta'(u^\epsilon)\,\lambda^\epsilon\,u^\epsilon+\eta'(u^\epsilon)\,\dv_{
g^\epsilon} f^\epsilon_\tau(u^\epsilon)=\epsilon\,\Delta_{g^\epsilon} \eta(u^\epsilon) -
\eta''(u_\epsilon)|du^\epsilon|_{g^\epsilon}\]
by a test function $\varphi$ as in Definition \ref{def:mvsol}, integrating over $M\times I$, and
integrating by parts we obtain
\begin{equation*}
 \begin{aligned}
&-\int_I\int_M \pa_t\varphi\,\eta(u^\epsilon)+ \varphi\,[\lambda^\epsilon\,(\eta(u^\epsilon)-\eta'(u^\epsilon)\,u^\epsilon) + \dv_{g^\epsilon}^x q^\epsilon_\tau(u^\epsilon)\\
&\hspace{1.4cm}-\eta'(u^\epsilon)\,\dv_{g^\epsilon}^x f^\epsilon_\tau(u^\epsilon)]
+ d\varphi(q^\epsilon_\tau(u^\epsilon))\ dV_{g^\epsilon}\,dt + \int_M\varphi(0)\,\eta(u_0^\epsilon)\
dV_{g^\epsilon(0)}\\
&\hspace{7.5cm}\le \int_I\int_M \epsilon\,\Delta_{g^\epsilon}\varphi\,\eta(u^\epsilon)\ dV_{g^\epsilon}\,dt.
 \end{aligned}
\end{equation*}
Letting $\epsilon\searrow 0$ we deduce that $\nu$ is a measure-valued entropy solution of \eqref{eqn:conlaw} with $f$ replaced by $f_\tau$. 

Now we let $\nu$ and $\sigma$ be two measure-valued entropy solutions for the same initial datum
$u_0\in L^\infty(M)$. By Proposition \ref{prop:acp}, for $\eta(u,v):=|u-v|$ and almost all
$0<s\le t<T$ we have
\[\int_M\langle\nu_t\otimes\sigma_t,\eta\rangle\ dV_{g(t)}\le
\int_M\langle\nu_s,|\,\cdot\,-u_0|\rangle + \langle\sigma_s,|\,\cdot\,-u_0|\rangle\ dV_{g(s)}.\] 
Lemma \ref{lem:timereg} shows that the left hand side vanishes for almost all $t\in I$, i.e. for
almost all $x,\,t$ the measures $\nu_{x,t}$ and $\sigma_{x,t}$ must concentrate at same point
$u(x,t)$. The resulting function $u$ is an entropy solution which is unique in the set of
measure-valued entropy solutions. 

By a standard argument from the theory of Young measures we see that $(u^\epsilon)$ converges to $u$ in $L^1(M\times I)$. In particular, a subsequence converges in $L^1(M)$ almost everywhere in $I$. Hence, in view of the weak-$^*$ lower semicontinuity of the $L^\infty(M)$ norm and the convergence of the constant $\norm{u^\epsilon_0}_{L^\infty(M)}$ and the functions $c_2=c_2^{\tau,\epsilon},c_3=c_3^{\epsilon}$ for $\epsilon\searrow 0$, we deduce that $u$ admits the bound \eqref{eqn:uniformbound} with $f$ replaced by $f_\tau$. Since $c_2^\tau\le c_2$ and by the choice of $N$ the function $u$, in fact, solves the original equation and admits the original bound \eqref{eqn:uniformbound}.
\qed\\

\begin{remark}
Taking into account the uniqueness of entropy solutions we see from the proof that the whole sequence $(u^\epsilon)$ of vanishing viscosity approximations converges to $u$ strongly in $L^p(M\times I)$ for all $1\le p<\infty$ and weakly-$^*$ in $L^\infty(M\times I)$. Furthermore, for another entropy solution $\tilde u$ we deduce from Proposition
\ref{prop:acp} and Lemma \ref{lem:timereg} the $L^1$ contraction property
\begin{equation*}
 \begin{aligned}
  \int_M |u(t)-\tilde u(t)|\ dV_{g(t)}\le \int_M |u(s)-\tilde u(s)|\ dV_{g(s)}
 \end{aligned}
\end{equation*}
for almost all $0<s\le t<T$ and for $s=0$ and almost all $t\in I$ as well as the comparison
principle
\[u\le \tilde u \text{ a.e. in } I\times M\text{ provided }u(0)\le \tilde u(0) \text{ a.e. in }M.\]
\end{remark}

\vspace{0.5cm}
\section{Total Variation Estimates}
\label{sec:tvestimates}
Throughout this section we assume that $g,f\in C^2_{\loc}$, and we assume $\pa_u\dv f$ to be uniformly bounded from below in $M\times I\times\setR$.\footnote{Again, it might be possible to lower the regularity assumptions (slightly).} Furthermore, we let the function spaces be normed w.r.t. the time-dependent metric.

\begin{definition} 
The total variation of a function $u\in L^1(M)$ is defined as
\begin{equation*}
 \tv_{g(t)}(u) := \sup \int_M u\, \dv_{g(t)} X \ dV_{g(t)}.
\end{equation*}
The supremum is taken over all smooth vector fields $X$ with
$\norm{X}_{L^\infty(M,g(t))}\le 1$.
%\footnote{Note that the norm of a vector field depends on the metric.}
For $u\in W^{1,1}(M)$ we have
\begin{equation*}
 \tv_{g(t)}(u) = \int_M |\nabla u|_{g(t)} \,dV_{g(t)}.
\end{equation*}
Furthermore, we let
 \begin{equation*}
  \BV(M) := \{u\in L^1(M)\ |\ \tv_{g(0)}(u) < \infty \}.
 \end{equation*}
\end{definition}
\vspace{1mm}
\begin{remark}\label{rem:BVwelldefined}
 The definition of the space $\BV(M)$ does not depend on the choice of the metric. To see this, let $g_1,g_2$ be two metrics of class $C^{0,1}$. Note that there is a unique, positive function $v\in C^{0,1}(M)$ such that
$dV_{g_1}=v\,dV_{g_2}$. Locally we have $v^2=\det ((g_1)_{ij})/\det ((g_2)_{ij})$. Hence
\begin{equation}\label{eq:tvequiv}
\begin{aligned}
 \int_M u\, \dv_{g_1} X \ dV_{g_1} =& \int_M u\, (\dv_{g_2} (v\,X)-dv(X)) \ dV_{g_2}\\
&+\int_M u\,(\dv_{g_1} - \dv_{g_2}) X\ dV_{g_1}.
\end{aligned}
\end{equation}
Since the difference $\dv_{g_1} - \dv_{g_2}$ is not acting as a derivative on $X$, the modulus of the right hand side of \eqref{eq:tvequiv} is dominated by
\begin{equation}\label{eq:tvequiv2}
\tv_{g_2}(u)\norm{v\,X}_{L^\infty(M,g_2)}+ c\norm{u}_{L^1(M,g_1)}\norm{X}_{L^\infty(M,g_1)}.
\end{equation}
This proves the claim.
% Then, their induced $L^\infty$-norms $\norm{\cdot}_{L^\infty(M,{g_1})}$, $\norm{\cdot}_{L^\infty(M,{g_2})}$
% are equivalent and the proportionality factor $\alpha$ of the volume forms
%  defined by $dV_{g_1}=\alpha \ dV_{g_2}$ is of class $C^{0,1}$.
% Furthermore, the difference between the respective divergences of a smooth vector field $X$ can be estimated via
% \begin{equation*}
%  \norm{\dv_{g_1}(X) - \dv_{g_2}(X)}_{L^\infty(M)} \leq  C \norm{X}_{L^\infty(M,g_i)}, \qquad {i=1,2}
% \end{equation*}
% for some constant $C$. Assuming $u\in L^1(M)$ with $\tv_{g_2}(u) < \infty$, then,
% for a smooth vector field $X$ with $\norm{X}_{L^\infty(M,g_1)}\leq 1$ and $\alpha_\infty:=\norm{\alpha}_{L^\infty(M)}$, 
% $X_{\infty,g^2}:=\norm{X}_{L^\infty(M,g_2)}$ we can estimate
% \begin{equation}\label{eq:tvequiv}
% \aligned
%  \int_M u \dv_{g_1} X \ dV_{g_1} =&
%  \alpha_\infty X_{\infty,g^2} \int_M u  \dv_{g_2}\left ( \frac{\alpha }{\alpha_\infty X_{\infty,g^2}} X \right) dV_{g_2} \\
% &- \int_M u\  g_2( \nabla_{g_2} \alpha, X ) \ dV_{g_2} \\
% &+ \int_M u  (\dv_{g_1} X - \dv_{g_2} X) \ \alpha \ dV_{g_2}\\
% &\leq C_1\ \tv_{g_2}(u) + C_2 \ \norm{u}_{L^1(M,g_2)} < \infty,
% \endaligned
% \end{equation}
% for constants $C_1$, $C_2$ that depend on the metrics $g_1$, $g_2$, but not on $X$.
% %Here, we used the fact that $\BV(M)$ functions permit, as in the Euclidean case, distributional derivatives (cf. \cite{MPP07}).
% Taking the supremum over all such vector fields $X$ completes the remark.
\end{remark}

We shall now derive estimates of the total variation of smooth solutions of the parabolic regularization. Again, it is crucial to work in a global, coordinate-independent way in order to get clean calculations and sharp estimates.
\begin{proposition}\label{prop:tvestreg}
Assume that $g,f,u_0$ are smooth and let $u$ be the smooth solution of \eqref{eqn:pareg} for
some fixed $\epsilon>0$.
Then we have for all $t\in\bar I$
\begin{equation}\label{eqn:tvestreg}
\begin{aligned}
 \tv_{g(t)}(u(t)) \le \Big(e^{\int_0^tc_4(s)\,ds}\,\tv_{g(0)}(u_0)+\int_0^t e^{\int_s^t c_4(\tau)\,d\tau}\,c_5(s)\ ds\Big)\\ \exp\big(\epsilon\norm{\ric}_{L^1(I,L^\infty(M))}\big),
\end{aligned}
\end{equation}
where
\begin{equation*}
\begin{aligned}
c_4(t)&:=-\inf_{x, X}\pa_t g_{x,t}(X,X) - \inf_{x,X,\bar u}\langle\nabla_X\,
\pa_u f (x,t,\bar u),X\rangle,\\
c_5(t)&:= u_{\max}\,\|\nabla\lambda(\,\cdot\,,t)\|_{L^1(M)} + \big\|\sup_{|\bar u|\le u_{\max}}|\nabla
\dv f(\,\cdot\,,t,\bar u)|\big\|_{L^1(M)}.
\end{aligned}
\end{equation*}
% \begin{equation*}
% \begin{aligned}
% c_4(t)&:= \tv_{g(0)}(u_0) + u_{\max}\,\|\nabla\lambda\|_{L^1(I\times M)} + \big\|\sup_{\bar u}|\nabla
% \dv f(\bar u)|\big\|_{L^1(I\times M)},\\
% c_5&:=\big\|\sup_{X}|\pa_t g(X,X)|\big\|_{L^1(I,L^\infty(M))} + \big\|\sup_{X,\,\bar u}|\langle\nabla_X\,
% \pa_u f (\bar u),X\rangle|\big\|_{L^1(I,L^\infty(M))}.
% % c_4&:= \tv(0,u_0) + \big\|u_{\max}\,|\nabla\lambda| + \max_{|\bar u|\le u_{\max}}|\nabla
% % \dv f|(\bar u)\big\|_{L^1(I\times M)},\\
% % c_5^\epsilon&:=\big\|\max_{|X|\le 1}\frac{|\pa_t g(X,X)|}{2}+ \max_{|X|\le 1,\,|\bar u|\le
% % u_{\max}}|\langle\nabla_X\,
% % \pa_u f(\bar u),X\rangle| + \epsilon\,|\ric|\big\|_{L^1(I,L^\infty(M))}.
% \end{aligned}
% \end{equation*}
Here $u_{\max}$ denotes maximum of the right hand side of \eqref{eqn:uniformbound} and the infima are taken over all $x\in M$, all tangent vectors $|X|_{g(t)}\le 1$, and (for the second infimum) all real numbers $|\bar u|\le u_{\max}$.
\end{proposition}
\proof 
Let us assume for the moment $g$ to be independent of time, in particular $\lambda\equiv 0$. Taking
the total covariant derivative of \eqref{eqn:pareg} we obtain
\begin{equation}\label{eqn:eqdiff}
\begin{aligned}
\pa_t\nabla u + \nabla\dv f = \epsilon\, \nabla \Delta u = \epsilon\,\big(\Delta\nabla u -
[\Delta,\nabla]u\big)=\epsilon\,\big(\Delta\nabla u - \ric(\nabla u,\,\cdot\,)\big).
\end{aligned}
\end{equation}
Here $\ric$ denotes the Ricci tensor. The commutator identity that is used in the last equality
follows from the definition of the Riemannian curvature tensor $R$ and
\[u^{\; i}_{;\; ij}-u^{\; i}_{;\; ji}=R_{ijk}^{\ \ \ i}\,u_;^{\, k}=\ric_{jk}\,u_;^{\,
k}.\]
In the following we will use the function $\eta_\delta$ from the proof of Proposition
\ref{prop:uniformbound}. Taking the scalar product of \eqref{eqn:eqdiff} with $\frac{\nabla
u}{|\nabla u|}\, \eta_\delta'(|\nabla u|)$ gives
\begin{equation}\label{eqn:zw}
\begin{aligned}
\pa_t \eta_\delta(|\nabla u|) + \langle\nabla\dv f(u),&\nabla u\rangle\,\frac{\eta_\delta'(|\nabla
u|)}{|\nabla u|}\\
&=\epsilon\,\frac{\eta_\delta'(|\nabla u|)}{|\nabla u|}\,\big(\langle\Delta\nabla u,\nabla u\rangle
- \ric(\nabla u,\nabla u)\big).
\end{aligned}
\end{equation}
Let us manipulate the second term on the left hand side. We have
\begin{equation*}
 \left(\dv f(u)\right)_{;i} = (f^j(u))_{;ji} =  \pa_u f^j_{\ ;i}(u)\,u_{;j} + f^j_{\ ;ji}(u) + 
 (\pa_u f^j(u)\, u_{;i})_{;j}.
\end{equation*}
Concerning the first term on the right hand side of this equation we remark that
\begin{equation*}
\begin{aligned}
 \pa_u f^j_{\ ;i}(u)\,u_{;j}\,u_{;}^{\, i}\ \frac{\eta_\delta'(|\nabla u|)}{|\nabla
u|}=\big\langle\big(\nabla_\frac{\grad u}{|\grad u|}\pa_u f\big)(u),\grad u\big\rangle\
\eta_\delta'(|\nabla u|),
\end{aligned}
\end{equation*}
while with regard to the last term we compute
\begin{equation*}
\begin{aligned}
 (\pa_u f^j(u)\,u_{;i})_{;j}\, u^{\, i}_{;}\ \frac{\eta_\delta'(|\nabla u|)}{|\nabla u|}
&= \dv \big(\pa_uf(u)\, \eta_\delta(|\nabla u|)\big) \\
&\hspace{0.5cm}+ \dv (\pa_u f(u))\, (|\nabla u| \, \eta_\delta'(|\nabla u|) - 
\eta_\delta(|\nabla u|)).
\end{aligned}
\end{equation*}
% Let us now analyze the first term on the right hand side of \eqref{eqn:zw}. Integrating over $M$
% and integrating by parts gives
% \begin{equation*}
% \begin{aligned}
% \epsilon\,\int_M\frac{\eta_\delta'(|\nabla u|)}{|\nabla u|}\,\big(\langle\Delta\nabla u,\nabla
% u\rangle % - \ric(\nabla u,\nabla u)\big).
% \end{aligned}
% \end{equation*}
Hence, integrating \eqref{eqn:zw} over $M$ results in
\begin{equation}\label{eqn:zw2}
\begin{aligned}
\frac{d}{dt}\int_M \eta_\delta(|\nabla u|)\ dV + \int_M &\big\langle\big(\nabla_\frac{\grad
u}{|\grad u|}\pa_u f\big)(u),\grad u\big\rangle\ \eta_\delta'(|\nabla u|)\\
&+ \big(\nabla_{\frac{\grad u}{|\grad u|}}\dv f\big)(u)\ \eta_\delta'(|\nabla u|) \\
&+ \dv (\pa_u f(u))\,(|\nabla u| \, \eta_\delta'(|\nabla u|) - 
\eta_\delta(|\nabla u|))\ dV\\
&\hspace{-0.8cm}=\epsilon\,\int_M\frac{\eta_\delta'(|\nabla u|)}{|\nabla
u|}\,\big(\langle\Delta\nabla u,\nabla u\rangle - \ric(\nabla u,\nabla u)\big)\ dV.
\end{aligned}
\end{equation}
Integration by parts (for general tensor fields, see, e.g., \cite{Lee97}) yields
\begin{equation*}
\int_M\frac{\eta_\delta'(|\nabla u|)}{|\nabla u|}\,\langle\Delta\nabla u,\nabla
u\rangle\ dV = - \int_M\big\langle\nabla^2 u,\nabla\Big(\frac{\eta_\delta'(|\nabla u|)}{|\nabla u|}
\nabla u\Big)\big\rangle\ dV.
\end{equation*}
Let us now show that the integrand on the right hand side is nonnegative. To this end we let
$G:\setR^d\rightarrow\setR$ denote the convex function $x\mapsto
\eta_\delta(|x|_{\setR^d})$, where $|x|_{\setR^d}$ is the standard Euclidean norm. Choosing
Riemannian normal coordinates centred at some point $p\in M$, the integrand evaluated in $p$ equals
\begin{equation*}
\begin{aligned}
\sum_{i,j=1}^d \pa_i\pa_j u \, \pa_j(\pa_i G)(\nabla u)=\sum_{i,j=1}^d\pa_i\pa_j u \, \pa_k\pa_j u\,
(\pa_k \pa_i G)(\nabla u)\geq 0.
\end{aligned}
\end{equation*}
Consequently, letting $\delta\searrow 0$ in \eqref{eqn:zw2} we obtain
\begin{equation}\label{eqn:zw3}
\begin{aligned}
\frac{d}{dt}\int_M |\nabla u|\ dV + \int_M \big\langle\big(\nabla_\frac{\grad u}{|\grad
u|}\pa_u f\big)(u),\grad u\big\rangle + \big(\nabla_{\frac{\grad u}{|\grad u|}}\dv f\big)(u)\ dV\\
\le-\epsilon\,\int_M\ric\big(\nabla u,\frac{\nabla u}{|\nabla u|}\big)\ dV
\end{aligned}
\end{equation}
almost everywhere in $I$ with $\frac{\nabla u}{|\nabla u|}:=0$ for $\nabla u=0$. An application of Gronwall's Lemma completes the proof in the case of a metric independent of time.

Now, let us drop this assumption on $g$. Then, instead of \eqref{eqn:zw}
we have
\begin{equation}\label{eqn:zw4}
\begin{aligned}
\langle\pa_t \nabla u,\nabla u\rangle\,\frac{\eta_\delta'(|\nabla
u|)}{|\nabla u|} + \langle \nabla (u\,\lambda),\nabla u\rangle\,\frac{\eta_\delta'(|\nabla
u|)}{|\nabla u|} +
\langle\nabla\dv f(u),\nabla u\rangle\frac{\eta_\delta'(|\nabla
u|)}{|\nabla u|}\\
=\epsilon\,\frac{\eta_\delta'(|\nabla u|)}{|\nabla
u|}\,\big(\langle\Delta\nabla u,\nabla u\rangle
- \ric(\nabla u,\nabla u)\big).
\end{aligned}
\end{equation}
The last term on the left hand side and the right hand side may be treated exactly like in
the time-independent setting. The only difference stems from the first two terms on the left hand
side. But these two terms equal
\begin{equation*}
\begin{aligned}
\pa_t \eta_\delta(|\nabla u|) + (\pa_t g)\big(\grad u,\frac{\grad
u}{2|\grad u|}\big)\,\eta_\delta'(|\nabla
u|) + \lambda\,|\nabla u|\,\eta_\delta'(|\nabla u|)\\
 + u\,\big\langle\nabla\lambda,\frac{\nabla
u}{|\nabla u|}\big\rangle\,\eta_\delta'(|\nabla u|)\\
=\pa_t \eta_\delta(|\nabla u|) + \lambda\,\eta_\delta(|\nabla u|) + (\pa_t g)\big(\grad u,\frac{\grad
u}{2|\grad u|}\big)\,\eta_\delta'(|\nabla u|)\\
 + u\,\big\langle\nabla\lambda,\frac{\nabla u}{|\nabla u|}\big\rangle\,\eta_\delta'(|\nabla u|) + \lambda\,\big(|\nabla
u|\,\eta_\delta'(|\nabla u|)-\eta_\delta(|\nabla u|)\big).
\end{aligned}
\end{equation*}
The last term tends to $0$ pointwise. Hence, integrating \eqref{eqn:zw4} over $M$ and letting $\delta\searrow
0$, instead of \eqref{eqn:zw3} we obtain
\begin{equation*}
\begin{aligned}
&\frac{d}{dt}\int_M |\nabla u|\ dV + \int_M \big\langle\big(\nabla_\frac{\grad u}{|\grad
u|}\pa_u f\big)(u),\grad u\big\rangle + \big(\nabla_{\frac{\grad u}{|\grad u|}}\dv f\big)(u)\\
&+ (\pa_t g)\big(\grad u,\frac{\grad
u}{2|\grad u|}\big) + u\,\big\langle\nabla\lambda,\frac{\nabla
u}{|\nabla u|}\big\rangle\ dV \le-\epsilon\,\int_M\ric\big(\nabla u,\frac{\nabla u}{|\nabla
u|}\big)\ dV
\end{aligned}
\end{equation*}
almost everywhere in $I$. Again, an application of Gronwall's Lemma proves the claim. \qed\\

Concerning the second part of the following theorem, we let $v\in C^{1,1}(M\times \bar I)$ denote the unique, positive function such that $dV_{g(t)}=v(t)\,dV_{g(0)}$.

\begin{theorem}\label{thm:tvest}
For every $u_0\in L^\infty(M)\cap\BV(M)$ the unique entropy solution $u\in L^\infty(M\times I)$
of \eqref{eqn:conlaw} satisfies for almost all $t\in I$
\begin{equation}\label{eqn:tv}
\tv_{g(t)}(u(t)) \le e^{\int_0^tc_4(s)\,ds}\,\tv_{g(0)}(u_0)+\int_0^t c_5(s)\,e^{\int_s^t c_4(\tau)\,d\tau}\ ds.
\end{equation}
% where $c_4,\,u_{\max}$ are taken from Proposition \ref{prop:tvestreg} and
% \[c_5:=\big\|\max_{|X|\le 1}\frac{|\pa_t g(X,X)|}{2}+ \max_{|X|\le 1,\,|\bar u|\le
% u_{\max}}|\langle\nabla_X\,
% \pa_u f(\bar u),X\rangle|\big\|_{L^1(I,L^\infty(M))}.\]
Furthermore, there exists a representative of $u$ such that for all $s,t\in \bar I$
\[\norm{u(t)\,v(t)-u(s)\,v(s)}_{L^1(M,g(0))} \le (c_6 + c_7\,\tv_{\max})\,|t-s|,\]
where $\tv_{\max}$ denotes the maximum of the right hand side of \eqref{eqn:tv} and
\[c_6:=\big\|\sup_{|\bar u|\le u_{\max}}|\dv f(\bar u)|\big\|_{L^\infty(I,L^1(M))},\ c_7:= \big\|\sup_{|\bar u|\le u_{\max}}|\pa_u
f|\big\|_{L^\infty(I\times M)}.\]
% \[c_6:=\big\|\max_{|\bar u|\le
% u_{\max}}|\dv f|(\bar u)\big\|_{L^1(I\times M)}
% + \big\|\max_{|\bar u|\le u_{\max}}|\pa_u
% f|(\bar u)\big\|_{L^1(I,L^\infty(M))}.\]
In particular, we have $u\in C^{0,1}(\bar I,L^1(M))$.  
\end{theorem}
\proof Using the procedure introduced in Section \ref{sec:not} we define mollifications $g^{\epsilon_0},f^{\epsilon_0}$ of $g,f$. One needs to proceed more carefully in order to produce mollifications of $u_0$. Let us fix some $\epsilon_0>0$. We let $u_0^{\epsilon_1}$ denote the evaluation at time $\epsilon_1$ of the heat flow w.r.t. the metric $g^{\epsilon_0}(0)$ and the initial value $u_0$. In particular, we have $\norm{u_0^{\epsilon_1}}_{L^\infty}\le \norm{u_0}_{L^\infty}$, and the sequence $(u_0^{\epsilon_1})$ tends to $u_0$ in $L^1(M)$. Furthermore, in \cite{MPP07} it is shown that the sequence $(\tv_{g^{\epsilon_0}(0)}(u_0^{\epsilon_1}))_{\epsilon_1}$ converges to $\tv_{g^{\epsilon_0}(0)}(u_0)$.
% mollifications of $u_0$ constructed in \cite[Proposition 1.4]{MPP07} w.r.t. the metric $g^{\epsilon_0}$. In particular, the sequence $(u_0^{\epsilon_1})$ converges to $u_0$ in $L^1(M)$ and the sequence $(\tv_{g^{\epsilon_0}}(0,u_0^{\epsilon_1}))_{\epsilon_1}$ tends to $\tv_{g^{\epsilon_0}}(0,u_0)$. 
By Proposition \ref{prop:smoothsol} there exist smooth solutions $u^{\epsilon}_{\epsilon_0,\epsilon_1}$ of \eqref{eqn:pareg} admitting the bound \eqref{eqn:tvestreg}, with the original data replaced by the mollifications. We let first $\epsilon_1$, then $\epsilon$, and finally $\epsilon_0$ tend to zero, the distinction between $\epsilon$ and $\epsilon_0$ being necessary in order to avoid a blow-up of the Ricci tensor on the right hand side of \eqref{eqn:tvestreg}. So let us start with the first limiting process. Inequality \eqref{eqn:tvestreg} shows that the solutions are uniformly bounded in $L^\infty(I,W^{1,1}(M))$, and \eqref{eqn:pareg} shows that they admit distributional time-derivatives which are uniformly bounded in $L^\infty(I,(W^{2,\infty}(M))')$. Since 
\[W^{1,1}(M)\embedding\embedding L^1(M)\embedding (W^{2,\infty}(M))',\]
by the classical Aubin-Lions theorem, see, e.g., \cite{Aub63}, a subsequence of $(u^{\epsilon}_{\epsilon_0,\epsilon_1})_{\epsilon_1}$ converges pointwise almost everywhere to a weak solution $u^{\epsilon}_{\epsilon_0}$ of \eqref{eqn:pareg} with $g,f$ replaced by $g^{\epsilon_0},f^{\epsilon_0}$. Taking into account the lower semicontinuity of the total variation we easily show that the estimate \eqref{eqn:tvestreg} with $g,f$ replaced by $g^{\epsilon_0},f^{\epsilon_0}$ holds for $u^{\epsilon}_{\epsilon_0}$. By another two, very similar applications of the Aubin-Lions theorem we let $\epsilon$ and $\epsilon_0$ tend to zero. In doing so we take into account the facts that
% the lower semicontinuity of the total variation
\begin{equation}\label{eqn:facts}
 \begin{aligned}
\tv_{g(t)}(u(t))&\le\liminf_{\epsilon_0}\tv_{g^{\epsilon_0}(t)}(u_{\epsilon_0}(t)),\\
\tv_{g^{\epsilon_0}(0)}(u_0)&\rightarrow \tv_{g(0)}(u_0)\text{ for }\epsilon_0\searrow 0,
 \end{aligned}
\end{equation}
and the convergence of the functions $c_4=c_4^{\epsilon_0},c_5=c_5^{\epsilon_0}$. Here, the convergence \eqref{eqn:facts}$_2$ follows from the argument employed in Remark \ref{rem:BVwelldefined}, since in the estimate analogous to \eqref{eq:tvequiv2} the functions $v=v^{\epsilon_0}$ tend to $1$ uniformly and the constant $c=c^{\epsilon_0}$ is of the form $\mathcal{O}(\epsilon_0)$.
% the convergence of $(\tv_{g^{\epsilon_0}(0)}(u_0))_{\epsilon_0}$ to $\tv_{g(0)}(u_0)$ which follows easily from the local identity
% \begin{equation*}
% \begin{aligned}
% \dv_{g^{\epsilon_0}(0)} X\,&dV_{g^{\epsilon_0}(0)}-\dv_{g(0)} X\,dV_{g(0)}\\
% &=X^k(r)\,\pa_k\,\big((\det(g^{\epsilon_0}_{ij}(r,0)))^{\frac 1 2}-(\det(g_{ij}(r,0)))^{\frac 1 2}\big)\,dr.
%  \end{aligned}
% \end{equation*}
We thus obtain the first part of the theorem.\footnote{Of course, instead of invoking the Aubin-Lions theorem we could as well proceed like in the last section to deduce strong convergence.}

In order to show the second part, let us use in \eqref{eqn:mvsol} with $\eta=\id$ the test
function \[\varphi=\varphi_{t,\epsilon}(\tau)\,\varphi_{s,\epsilon}(2s-\tau)\,\psi(x)\] with $\varphi_{t,\epsilon}$ and $\psi$ as
in the proof of Lemma \ref{lem:timereg}. Letting
$\epsilon\searrow 0$ we obtain for almost all $0< s\le t< T$
\begin{equation}\label{eqn:lipreg}
\begin{aligned}
\int_M (u(t)\,v(t)-u(s)\,v(s))\,\psi\ dV_{g(0)}&=\int_M u(t)\,\psi\
dV_{g(t)} - \int_M u(s)\,\psi\ dV_{g(s)}\\
& = - \int_s^t\int_M d\psi(f(u))\
dV\,d\tau.
\end{aligned}
\end{equation}
Another approximation argument (using again the result from \cite{MPP07}) and the chain rule show that the modulus of the right hand side is dominated by
\[\Big(c_6\, |t-s| + c_7\,\int_s^t \tv_{g(\tau)}(u(\tau))\ d\tau\Big)\norm{\psi}_{L^\infty(M)}.\]
Taking in \eqref{eqn:lipreg} the supremum over all $\norm{\psi}_{L^\infty(M)}\le 1$ completes the
proof.
\qed\\

\begin{remark}
If we assume that $g$ is independent of time, $\dv f(u=0)=0$ and 
\begin{equation}\label{eqn:symm}
\langle\nabla_X \pa_u f(x,t,\bar u),X\rangle=0
\end{equation}
for all $(x,t,\bar u)\in M\times I\times\setR$ and all
tangent vectors $X\in T_xM$, then
\eqref{eqn:tv} is a total variation diminishing estimate, i.e.
\begin{equation}\label{eqn:tvd}
\tv_g(u(t))\le\tv_g(u_0).
\end{equation}
Note that from \eqref{eqn:symm} we can deduce that $\pa_u\dv  f = \pa_u f_{\ ;i}^{i}=0$, in particular $\dv f=0$.
Condition \eqref{eqn:symm} says that the symmetric part of
the differential vanishes. Obviously this is equivalent to the claim that the differential
is skew-sym\-metric, i.e.
\begin{equation*}
\langle\nabla_X \pa_u f,Y\rangle=-\langle\nabla_Y \pa_u f,X\rangle
\end{equation*}
for all tangent vectors $X,\,Y$. Vector fields satisfying this condition are known in geometry
as Killing fields or infinitesimal isometries, because the flow generated by such a field is a one-parameter group of isometries.

Obviously, the assumption
\begin{equation*}
\langle\nabla_X \pa_u f,X\rangle\ge 0
\end{equation*}
would suffice to obtain \eqref{eqn:tvd}. But on closed manifolds this seemingly weaker condition is equivalent to
\eqref{eqn:symm}, since it implies
\[\dv \pa_u f = (\pa_u f)_{\; ;i}^{i} \ge 0.\]
By the Stokes theorem this expression must vanish. But then for each fixed $1\le i\le d$ we have $(\pa_u f)_{\; ;i}^{i}=0$, which
gives \eqref{eqn:symm}.
\end{remark}
\vspace{1mm}

\begin{remark}
If we assume $g$ to be time-independent, then, essentially, \eqref{eqn:symm} is not only sufficient, but also necessary for \eqref{eqn:tvd} to hold. Let us assume for simplicity that $g,f$ are smooth and that $\dv f=0$. If $\eqref{eqn:symm}$ does not hold, then we have $\langle\nabla_X \pa_u f(x,t,\bar u),X\rangle<0$ for some $x\in M,t\in\bar I,\bar u\in\setR,$ and some tangent vector $X\in T_xM$. We may assume without loss of generality that $t=0$. Let us fix a positively oriented, orthonormal basis of eigenvectors w.r.t. the symmetric part of the differential of $\pa_u f(0,x,\bar u)$, with the first eigenvector pointing in the direction of $X$. We denote by $\lambda_i$ the corresponding eigenvalues, in particular $\lambda_1<0$. This choice of basis induces normal coordinates $(r^i)$ centred at $x$. We set \[u_0(r):=\bar u + \kappa_0\,\rho(r^1/\kappa_1)\prod_{i=2}^d\rho(r^i/\kappa_2)\] for small constants $\kappa_0>0,\,\kappa_2>\kappa_1>0$ and the bump function $\rho$ introduced in Section \ref{sec:not}. 
Outside the coordinate patch we set $u_0\equiv \bar u$. For short 
times there exists a unique smooth solution $u$ for these initial values.\footnote{This can be easily shown by the usual implicit function theorem argument.} We shall choose the constants $\kappa_0,\kappa_1,\kappa_2$ such that the time-derivative
\begin{equation*}
\begin{aligned}
&\frac{d}{dt}\int_M |\nabla u|\ dV=-\int_M \big\langle\big(\nabla_\frac{\grad u}{|\grad u|}\pa_u f\big)(u),\grad u\big\rangle\ dV,
\end{aligned}
\end{equation*}
cf. the proof of Proposition \ref{prop:tvestreg}, is positive at $t=0$. Note that the right hand side is continuous in $u$ and can be written in the constructed normal coordinates at time $t=0$ as
\begin{equation}\label{eqn:ineq2}
\begin{aligned}
-\int_{\setR^d}\sum_{i=1}^d(\lambda_i+\mathcal{O}(\kappa_0)+\mathcal{O}(\kappa_2))\,\frac{|\pa_{r^i} u_0(r)|^2}{|\nabla u_0(r)|_{g(r,0)}}\,(1+\mathcal{O}(\kappa_2))\ dr.
\end{aligned}
\end{equation}
By the change of variables $\tilde r^1=r^1/\kappa_1$ and $\tilde r^i=r^i/\kappa_2$ for $i>1$ we easily see that there exists a constant $c>0$ such that
\[c\,\kappa_0\,\kappa_2^{d-1}\le\int_{\setR^d}\frac{|\pa_1 u_0|^2}{|\nabla
u_0|_{\setR^d}}\ dr\le c\,\kappa_0\,\kappa_2^d/\kappa_1\, \text{ and }\, \int_{\setR^d}\frac{|\pa_i u_0|^2}{|\nabla
u_0|_{\setR^d}}\ dr\le c\,\kappa_0\,\kappa_1\,\kappa_2^{d-2}\]
for $i>1$. Thus, taking into account that $g(r,0)$ is comparable to the Euclidean metric, we deduce that \eqref{eqn:ineq2} is positive if
\[1/c>(\kappa_0+\kappa_2)\, \kappa_2/\kappa_1 + (1+\kappa_0+\kappa_2)\,\kappa_1/\kappa_2\]
for some new, possibly large constant $c>0$. The right hand side can be chosen arbitrarily small by setting $\kappa_0=\kappa_2=\epsilon$ and $\kappa_1=\epsilon \sqrt{\epsilon}$ for $\epsilon >0$ small enough.

\end{remark}

\bibliographystyle{plain}    
\bibliography{literatur}

\end{document}